\documentclass[11pt,reqno]{amsart}
\usepackage{geometry}
\usepackage{graphicx}
\usepackage{amssymb}
\usepackage{caption}
\usepackage{subcaption}
\usepackage{enumitem}
\usepackage{mathtools}
\usepackage[all]{xy}
\usepackage{xcolor}
\usepackage{soul}
\usepackage[hidelinks,colorlinks=true,linkcolor=blue,urlcolor=blue,citecolor=blue]{hyperref}

\usepackage{microtype}

\usepackage{bbold}

\usepackage{MnSymbol}
\usepackage[linesnumbered,ruled,vlined]{algorithm2e}

\SetCommentSty{mycommfont}
\SetKwInput{KwInput}{Input}               
\SetKwInput{KwOutput}{Output}     

\usepackage{etoolbox} 

\makeatletter
\patchcmd{\mmeasure@}{\measuring@true}{
  \measuring@true
  \ifnum-\linenopenaltypar>\interdisplaylinepenalty
    \advance\interdisplaylinepenalty-\linenopenalty
  \fi
  }{}{}
\makeatother

\usepackage{graphicx}

\graphicspath{{figures/}}

\numberwithin{equation}{section}

\newtheorem{theorem}{Theorem}[section]

\newtheorem*{remark*}{Remark}

\usepackage[numbers,sort&compress]{natbib}
\bibpunct{[}{]}{;}{n}{,}{,}

\usepackage{dutchcal}

\DeclarePairedDelimiterX{\inner}[2]{\langle}{\rangle}{#1, #2}

\DeclareMathOperator*{\ext}{ext}

\allowdisplaybreaks

\title[Time-adaptive Lagrangian Variational Integrators for Accelerated Optimization]{Time-adaptive Lagrangian Variational Integrators\\ for Accelerated Optimization on Manifolds}
\author{Valentin Duruisseaux and Melvin Leok}

\makeatletter
\let\my@abstract=\relax
\def\abstract#1{%
  \def\my@abstract{%
    \normalfont\Small
    \list{}{\labelwidth\z@
      \leftmargin3pc \rightmargin\leftmargin
      \listparindent\normalparindent \itemindent\z@
      \parsep\z@ \@plus\p@
      \let\fullwidthdisplay\relax
    }%
    \item[\hskip\labelsep\scshape\abstractname.]%
    #1
  \endlist}}
\def\@setabstracta{%
  \ifx\my@abstract\relax
  \else
    \skip@20\p@ \advance\skip@-\lastskip
    \advance\skip@-\baselineskip \vskip\skip@
  \my@abstract
    \prevdepth\z@ % because \abstractbox is a vtop
  \fi
}
\makeatother

\begin{document}

\abstract{
A variational framework for accelerated optimization was recently introduced on normed vector spaces and Riemannian manifolds in~\citet{WiWiJo16} and~\citet{Duruisseaux2022Riemannian}. It was observed that a careful combination of time-adaptivity and symplecticity in the numerical integration can result in a significant gain in computational efficiency. It is however well known that symplectic integrators lose their near energy preservation properties when variable time-steps are used. The most common approach to circumvent this problem involves the Poincar\'e transformation on the Hamiltonian side, and was used in~\citet{duruisseaux2020adaptive} to construct efficient explicit algorithms for symplectic accelerated optimization. However, the current formulations of Hamiltonian variational integrators do not make intrinsic sense on more general spaces such as Riemannian manifolds and Lie groups. In contrast, Lagrangian variational integrators are well-defined on manifolds, so we develop here a framework for time-adaptivity in Lagrangian variational integrators and use the resulting geometric integrators to solve optimization problems on vector spaces and Lie groups.
}

\maketitle

\vspace{-7mm} 

\section{Introduction}

Many machine learning algorithms are designed around the minimization of a loss function or the maximization of a likelihood function. Due to the ever-growing scale of data sets, there has been a lot of focus on first-order optimization algorithms because of their low cost per iteration. In 1983, Nesterov's accelerated gradient method was introduced in ~\cite{Nes83}, and was shown to converge in $\mathcal{O}(1/k^2)$ to the minimum of the convex objective function $f$, improving on the $\mathcal{O}(1/k)$ convergence rate exhibited by the standard gradient descent methods.
This $\mathcal{O}(1/k^2)$ convergence rate was shown in ~\cite{Nes04} to be optimal among first-order methods using only information about $\nabla f$ at consecutive iterates. This phenomenon in which an algorithm displays this improved rate of convergence is referred to as acceleration, and other accelerated algorithms have been derived since Nesterov's algorithm. More recently, it was shown in ~\cite{SuBoCa16} that Nesterov's accelerated gradient method limits, as the time-step goes to 0, to a second-order differential equation and that the objective function $f(x(t))$ converges to its optimal value at a rate of $\mathcal{O}(1/t^2)$ along the trajectories of this ordinary differential equation. It was later shown in ~\cite{WiWiJo16} that in continuous time, the convergence rate of $f(x(t))$ can be accelerated to an arbitrary convergence rate $\mathcal{O}(1/t^p)$ in normed spaces, by considering flow maps generated by a family of time-dependent Bregman Lagrangian and Hamiltonian systems which is closed under time-rescaling. This framework for accelerated optimization in normed vector spaces has been studied and exploited using geometric numerical integrators in \cite{JordanSymplecticOptimization,Jordan2018,Campos2021,Franca2021,Franca2021b,duruisseaux2020adaptive}. In~\cite{duruisseaux2020adaptive}, time-adaptive geometric integrators have been proposed to take advantage of the time-rescaling property of the Bregman family and design efficient explicit algorithms for symplectic accelerated optimization. It was observed that a careful use of adaptivity and symplecticity could result in a significant gain in computational efficiency, by simulating higher-order Bregman dynamics using the computationally efficient lower-order Bregman integrators applied to the time-rescaled dynamics.

More generally, symplectic integrators form a class of geometric integrators of interest since, when applied to Hamiltonian systems, they yield discrete approximations of the flow that preserve the symplectic 2-form and as a result also preserve many qualitative aspects of the underlying dynamical system (see ~\cite{HaLuWa2006}). In particular, when applied to conservative Hamiltonian systems, symplectic integrators exhibit excellent long-time near-energy preservation. However, when symplectic integrators were first used in combination with variable time-steps, the near-energy preservation was lost and the integrators performed poorly (see ~\cite{CalSan93, GlaDunCan91}). There has been a great effort to circumvent this problem, and there have been many successes, including methods based on the Poincar\'e transformation~\cite{Zare1975,Ha1997}: a Poincar\'e transformed Hamiltonian in extended phase space is constructed which allows the use of variable time-steps in symplectic integrators without losing the nice conservation properties associated to these integrators. In~\cite{duruisseaux2020adaptive}, the Poincar\'e transformation was incorporated in the Hamiltonian variational integrator framework which provides a systematic method for constructing symplectic integrators of arbitrarily high-order based on the discretization of Hamilton's principle~\cite{MaWe2001, HaLe2012}, or equivalently, by the approximation of the generating function of the symplectic flow map. The Poincar\'e transformation was at the heart of the construction of time-adaptive geometric integrators for Bregman Hamiltonian systems which resulted in efficient, explicit algorithms for accelerated optimization in~\cite{duruisseaux2020adaptive}.

In~\cite{Tao2020,Lee2021}, accelerated optimization algorithms were proposed in the Lie group setting for specific choices of parameters in the Bregman family, and~\cite{alimisis2020} provided a first example of Bregman dynamics on Riemannian manifolds. The entire variational framework was later generalized to the Riemannian manifold setting in ~\cite{Duruisseaux2022Riemannian}, and time-adaptive geometric integrators taking advantage of the time-rescaling property of the Bregman family have been proposed in the Riemannian manifold setting as well using discrete variational integrators incorporating holonomic constraints~\cite{Duruisseaux2022Constrained} and projection-based variational integrators~\cite{Duruisseaux2022Projection}. Note that both these strategies relied on exploiting the structure of the Euclidean spaces in which the Riemannian manifolds are embedded. Although the Whitney and Nash Embedding Theorems~\cite{Whitney1944_2,Whitney1944_1,Nash1956} imply that there is no loss of generality when studying Riemannian manifolds only as submanifolds of Euclidean spaces, designing intrinsic methods that would exploit and preserve the symmetries and geometric properties of the manifold could have advantages both in terms of computational efficiency and in terms of improving our understanding of the acceleration phenomenon on Riemannian manifolds. Developing an intrinsic extension of Hamiltonian variational integrators to manifolds would require some additional work, since the current approach involves Type~II/III generating functions $H_d^+(q_k, p_{k+1})$, $H_d^-(p_k, q_{k+1})$, which depend on the position at one boundary point, and the momentum at the other boundary point. However, this does not make intrinsic sense on a manifold, since one needs the base point in order to specify the corresponding cotangent space. On the other hand, Lagrangian variational integrators involve a Type~I generating function $L_d(q_k, q_{k+1})$ which only depends on the position at the boundary points and is therefore well-defined on manifolds, and many Lagrangian variational integrators have been derived on Riemannian manifolds, especially in the Lie group~\cite{LeLeMc2007a,LeLeMc2007b, BRMa2009,HaLe2012,Nordkvist2010,Le2004,Hussein2006,LeLeMc2005,LeeThesis} and homogeneous space~\cite{LeLeMc2009b} settings. This gives an incentive to construct a mechanism on the Lagrangian side which mimics the Poincar\'e transformation, since it is more natural and easier to work on the Lagrangian side on more general spaces than on the Hamiltonian side. However, a first difficulty is that the Poincar\'e transformed Hamiltonian is degenerate and therefore does not have a corresponding Type~I Lagrangian formulation. As a result, we cannot exploit the usual correspondence between Hamiltonian and Lagrangian dynamics and need to come up with a different strategy. A second difficulty is that all the literature to this day on the Poincar\'e transformation constructs the Poincar\'e transformed system by reverse-engineering, which does not provide much insight into the origin of the mechanism and how it can be extended to different systems. \\

\noindent \textbf{Outline.} We first review the basics of variational integration of Lagrangian and Hamiltonian systems, and the Poincar\'e transformation in Section~\ref{section: Background}. We then introduce a simple but novel derivation of the Poincar\'e transformation from a variational principle in Section~\ref{section: Variational Derivation of Poincare Hamiltonian}. This gives additional insight into the transformation mechanism and provides natural candidates for time-adaptivity on the Lagrangian side, which we then construct both in continuous and  discrete time in Sections~\ref{section: Variational Time-Adaptivity Lagrangian} and \ref{section: Guess Time-Adaptivity Lagrangian}. We then compare the performance of the resulting time-adaptive Lagrangian accelerated optimization algorithms to their Poincar\'e Hamiltonian analogues in Section~\ref{section: Accelerated Optimization Vector Spaces}. Finally, we demonstrate in Section~ \ref{section: Accelerated Optimization General Spaces} that our time-adaptive Lagrangian approach extends naturally to more general spaces without having to face the obstructions experienced on the Hamiltonian side. \\

\noindent \textbf{Contributions.} In summary, the main contributions of this paper are:
\begin{itemize}
	\item A novel derivation of the Poincar\'e transformation from a variational principle, in Section~\ref{section: Variational Derivation of Poincare Hamiltonian}
	\item New frameworks for variable time-stepping in Lagrangian integrators, and new discrete variational formulations of Lagrangian mechanics with these new variable time-stepping mechanisms, in Sections~\ref{section: Variational Time-Adaptivity Lagrangian} and~\ref{section: Guess Time-Adaptivity Lagrangian}
	\item New explicit symplectic accelerated optimization algorithms on normed vector spaces 
	\item  New intrinsic symplectic accelerated optimization algorithms on Riemannian manifolds 
\end{itemize}

\hfill

\section{Background} \label{section: Background}

\subsection{Lagrangian and Hamiltonian Mechanics}  \label{Section: Lagrangian and Hamiltonian Mechanics}

Given a $n$-dimensional manifold $\mathcal{Q}$, a Lagrangian is a function $L:T\mathcal{Q}  \rightarrow \mathbb{R}$. The corresponding action integral $\mathcal{S}$ is the functional
\begin{equation} 
	\mathcal{S} (q) = \int_{0}^{T}{L(q,\dot{q})dt},
\end{equation}
over the space of smooth curves $q:[0,T] \rightarrow \mathcal{Q}$. Hamilton's variational principle states that $ \delta \mathcal{S}=0$ where the variation $\delta \mathcal{S}$ is induced by an infinitesimal variation $\delta q$ of the trajectory $q$ that vanishes at the endpoints. Given local coordinates $(q^1, \ldots , q^n)$ on the manifold $\mathcal{Q}$, Hamilton's variational principle can be shown to be equivalent to the Euler--Lagrange equations,
\begin{equation}\label{eq: EL Basic}
	\frac{d}{dt} \left( \frac{\partial L} {\partial \dot{q}^k} \right)= \frac{\partial L }{\partial q^k}, \qquad \text{for } k=1,\ldots , n. 
\end{equation}

The Legendre transform $\mathbb{F}L : T\mathcal{Q} \rightarrow T^* \mathcal{Q}$ of $L$ is defined fiberwise by $\mathbb{F}L : (q^i,\dot{q}^i) \mapsto \left(q^i,\frac{\partial L}{\partial \dot{q}^i} \right)$, and we say that a Lagrangian $L$ is regular or nondegenerate if the Hessian matrix $\frac{\partial^2 L}{\partial \dot{q}^2}$ is invertible for every $q$ and $\dot{q}$, and hyperregular if the Legendre transform $\mathbb{F}L$ is a diffeomorphism. A hyperregular Lagrangian on $T\mathcal{Q}$ induces a Hamiltonian system on $T^* \mathcal{Q}$ via 
\begin{equation}
	H(q,p) = \langle \mathbb{F}L (q,\dot{q}) , \dot{q} \rangle - L(q,\dot{q})= \sum_{j=1}^{n} {p_j  \dot{q}^j} - L(q,\dot{q}) \bigg|_{p_i=\frac{\partial L}{\partial \dot{q}^i}},
\end{equation} 
where $p_i = \frac{\partial L}{\partial \dot{q}^i} \in T^* \mathcal{Q}$ is the conjugate momentum of $q^i$. A Hamiltonian $H$ is called hyperregular if $\mathbb{F}H : T^* \mathcal{Q} \rightarrow T\mathcal{Q}$ defined by $
	\mathbb{F}H(\alpha) \cdot \beta = \frac{d}{ds} \big|_{s=0} H(\alpha +s\beta ) $, is a diffeomorphism. Hyperregularity of the Hamiltonian $H$ implies invertibility of the Hessian matrix $\frac{\partial^2 H}{\partial p^2}$ and thus nondegeneracy of~$H$. Theorem~7.4.3 in~\cite{MaRa1999} states that hyperregular Lagrangians and hyperregular Hamiltonians correspond in a bijective manner. We can also define a Hamiltonian variational principle on the Hamiltonian side in momentum phase space which is equivalent to Hamilton's equations, 
\begin{equation} \label{eq: Hamilton Equations Basic} \dot{p}_k = -\frac{\partial H}{\partial q^k} (p,q),  \qquad  \dot{q}^k = \frac{\partial H}{\partial p_k} (p,q), \qquad \text{for } k=1,\ldots , n.
\end{equation}
These equations can also be shown to be equivalent to the Euler--Lagrange equations \eqref{eq: EL Basic}, provided that the Lagrangian is hyperregular.

\subsection{Variational Integrators} \label{section: Variational Integrators}

Variational integrators are derived by discretizing Hamilton's principle, instead of discretizing Hamilton's equations directly. As a result, variational integrators are symplectic, preserve many invariants and momentum maps, and have excellent long-time near-energy preservation (see~\cite{MaWe2001}). 

Traditionally, variational integrators have been designed based on the Type~I generating function commonly known as the discrete Lagrangian, $L_d:Q \times Q \rightarrow \mathbb{R}$. The exact discrete Lagrangian that generates the time-$h$ flow of Hamilton's equations can be represented both in a variational form and in a boundary-value form. The latter is given by
\begin{equation}
	L_d^E(q_0,q_1;h)=\int_0^h L(q(t),\dot q(t)) dt  \label{exact_Ld},
\end{equation}
where $q(0)=q_0,$ $q(h)=q_1,$ and $q$ satisfies the Euler--Lagrange equations over the time interval $[0,h]$. A variational integrator is defined by constructing an approximation $L_d:Q \times Q \rightarrow \mathbb{R}$ to $L_d^E$, and then applying the discrete Euler--Lagrange equations,
\begin{equation}
	p_k=-D_1 L_d(q_k, q_{k+1}),\qquad p_{k+1}=D_2 L_d(q_k, q_{k+1}),  \label{IDEL}
\end{equation}
where $D_i$ denotes a partial derivative with respect to the $i$-th argument, and these equations implicitly define the integrator map $\tilde{F}_{L_d}:(q_k,p_k)\mapsto(q_{k+1},p_{k+1})$. The error analysis is greatly simplified via Theorem~2.3.1 of~\cite{MaWe2001}, which states that if a discrete Lagrangian, $L_d:Q\times Q\rightarrow\mathbb{R}$, approximates the exact discrete Lagrangian $L_d^E:Q\times Q\rightarrow\mathbb{R}$ to order $r$, i.e.,
\begin{equation} 
	L_d(q_0, q_1;h)=L_d^E(q_0,q_1;h)+\mathcal{O}(h^{r+1}) ,
\end{equation}
then the discrete Hamiltonian map $\tilde{F}_{L_d}:(q_k,p_k)\mapsto(q_{k+1},p_{k+1})$, viewed as a one-step method, has order of accuracy $r$. Many other properties of the integrator, such as momentum conservation properties of the method, can be determined by analyzing the associated discrete Lagrangian, as opposed to analyzing the integrator directly. 

Variational integrators have been extended to the framework of Type~II and Type~III generating functions, commonly referred to as discrete Hamiltonians (see~\cite{LaWe2006, LeZh2011,ScLe2017}). Hamiltonian variational integrators are derived by discretizing Hamilton's phase space principle. The boundary-value formulation of the exact Type~II generating function of the time-$h$ flow of Hamilton's equations is given by the exact discrete right Hamiltonian,
\begin{equation}
	H_d^{+,E}(q_0,p_1;h) =  p_1^\top q_1 - \int_0^h \left[ p(t)^\top \dot{q}(t)-H(q(t), p(t)) \right] dt, \label{exact_Hd}
\end{equation}
where $(q,p)$ satisfies Hamilton's equations with boundary conditions $q(0)=q_0$ and $p(h)=p_1$. A Type~II Hamiltonian variational integrator is constructed by using a discrete right Hamiltonian $H_d^+$ approximating $	H_d^{+,E}$, and applying the discrete right Hamilton's equations,
\begin{equation}\label{Discrete Right Eq}
	p_0=D_1H_d^+(q_0,p_1), \qquad q_1=D_2H_d^+(q_0,p_1),
\end{equation}
which implicitly defines the integrator, $\tilde{F}_{H_d^+}:(q_0,p_0) \mapsto (q_1,p_1)$.

Theorem~2.3.1 of~\cite{MaWe2001}, which simplified the error analysis for Lagrangian variational integrators, has an analogue for Hamiltonian variational integrators. Theorem~2.2 in~\cite{ScLe2017} states that if a discrete right Hamiltonian $H^+_d$ approximates the exact discrete right Hamiltonian $H_d^{+,E}$ to order $r$, i.e.,
\begin{equation} 
	H^+_d(q_0, p_1;h)=H_d^{+,E}(q_0,p_1;h)+\mathcal{O}(h^{r+1}),
\end{equation}
then the discrete right Hamilton's map $\tilde{F}_{H^+_d}:(q_k,p_k)\mapsto(q_{k+1},p_{k+1})$, viewed as a one-step method, is order $r$ accurate. Note that discrete left Hamiltonians and corresponding discrete left Hamilton's maps can also be constructed in the Type~III case (see~\cite{LeZh2011,duruisseaux2020adaptive}).

Examples of variational integrators include Galerkin variational integrators ~\cite{LeZh2011,MaWe2001}, Prolongation-Collocation variational integrators~\cite{LeSh2011}, and Taylor variational integrators~\cite{ScShLe2017}. In many cases, the Type~I and Type~II/III approaches will produce equivalent integrators. This equivalence has been established in~\cite{ScShLe2017} for Taylor variational integrators provided the Lagrangian is hyperregular, and in~\cite{LeZh2011} for generalized Galerkin variational integrators constructed using the same choices of basis functions and numerical quadrature formula provided the  Hamiltonian is hyperregular. However, Hamiltonian and Lagrangian variational integrators are not always equivalent. In particular, it was shown in~\cite{ScLe2017} that even when the Hamiltonian and Lagrangian integrators are analytically equivalent, they might still have different numerical properties because of numerical conditioning issues. Even more to the point, Lagrangian variational integrators cannot always be constructed when the underlying Hamiltonian is degenerate. This is particularly relevant in variational accelerated optimization since the time-adaptive Hamiltonian framework for accelerated optimization presented in~\cite{duruisseaux2020adaptive} relies on a degenerate Hamiltonian which has no associated Lagrangian description. We will thus not be able to exploit the usual correspondence between Hamiltonian and Lagrangian dynamics and will have to come up with a different strategy to allow time-adaptivity on the Lagrangian side.

We now describe the construction of Taylor variational integrators as introduced in~\cite{ScShLe2017} as we will use them in our numerical experiments. A discrete approximate Lagrangian or Hamiltonian is constructed by approximating the flow map and the trajectory associated with the boundary values using a Taylor method, and approximating the integral by a quadrature rule. The Taylor variational integrator is generated by the implicit discrete Euler--Lagrange equations associated with the discrete Lagrangian or by the Hamilton's equations associated with the discrete Hamiltonian. More explicitly, we first construct $\rho$-order and $(\rho+1)$-order Taylor methods $\Psi_h^{(\rho)}$ and $\Psi_h^{(\rho+1)}$ approximating the exact time-$h$ flow map $\Phi _h : TQ \rightarrow TQ$ corresponding to the Euler--Lagrange equation in the Type I case or the exact time-$h$ flow map $\Phi _h : T^*Q \rightarrow T^*Q$ corresponding to Hamilton's equation in the Type II case. Let $\pi_{Q}:(q,p)\mapsto q$ and $\pi_{T^*Q}:(q,p)\mapsto p$. Given a quadrature rule of order $s$ with weights and nodes $(b_i,c_i)$ for $i=1,...,m$, the Taylor variational integrators are then constructed as follows: \\

\noindent \underline{\textbf{Type I Lagrangian Taylor Variational Integrator (LTVI)}:}
\begin{enumerate}[label=(\roman*)]
	
	\item Approximate $\dot{q}(0)=v_0$ by the solution $\tilde{v}_0$ of the problem 
	$ q_1 = \pi_{Q} \circ \Psi_h^{(\rho+1)}(q_0,\tilde{v}_0) .$
	
	\item Generate approximations $(q_{c_i},v_{c_i}) \approx (q(c_i h),\dot{q}(c_i h))$ via
	$ (q_{c_i},v_{c_i})  = \Psi_{c_i h}^{(\rho)}(q_0,\tilde{v}_0).$
	
	\item Apply the quadrature rule to obtain the associated discrete Lagrangian,
	\begin{equation*}
		L_d(q_0,q_1;h) = h \sum_{i=1}^{m}{b_i L(q_{c_i} , v_{c_i})}.
	\end{equation*} 
	
	\item The variational integrator is then defined by the implicit discrete Euler--Lagrange equations, 
	\begin{equation*}
		p_0=-D_1 L_d(q_0, q_1),\qquad p_1=D_2 L_d(q_0, q_1).
	\end{equation*}
\end{enumerate} 
\hfill

\noindent \underline{\textbf{Type II Hamiltonian Taylor Variational Integrator (HTVI)}:}

\begin{enumerate}[label=(\roman*)]
	
	\item Approximate $p(0)=p_0$ by the solution $\tilde{p}_0$ of the problem 
	$ p_1 = \pi_{T^*Q} \circ \Psi_h^{(\rho)}(q_0,\tilde{p}_0) .$
	
	\item Generate approximations $(q_{c_i},p_{c_i}) \approx (q(c_i h),p(c_i h))$ via
	$ (q_{c_i},p_{c_i})  = \Psi_{c_i h}^{(\rho)}(q_0,\tilde{p}_0).$
	
	\item Approximate $q_1$ via $ \tilde{q}_1 = \pi_{Q}  \circ \Psi_h^{(\rho+1)}(q_0,\tilde{p}_0).$
	
	\item Use the continuous Legendre transform to obtain $\dot{q}_{c_i} = \frac{\partial H}{\partial p_{c_i}}$.
	
	\item Apply the quadrature rule to obtain the associated discrete right Hamiltonian,
	\begin{equation*}
		H_d^+(q_0,p_1;h) = p_1^{\top} \tilde{q}_1 - h \sum_{i=1}^{m}{b_i \left[  p_{c_i}^{\top} \dot{q}_{c_i} - H(q_{c_i},p_{c_i})  \right]}.
	\end{equation*} 
	
	\item The variational integrator is then defined by the implicit discrete right Hamilton's equations, 
	\[ q_1 = D_2 H_d^+(q_0,p_1), \qquad p_0 = D_1 H_d^+(q_0,p_1).\] 
\end{enumerate} 
\hfill 

\noindent The following error analysis results were derived in~\cite{ScShLe2017} and~\cite{duruisseaux2020adaptive}:
\begin{theorem}
	Suppose the Lagrangian $L$ is Lipschitz continuous in both variables, and is sufficiently regular for the Taylor method $\Psi_h^{(\rho+1)}$ to be well-defined. 
	
	\noindent Then $L_d(q_0,q_1)$ approximates $L_d^{E}(q_0,q_1)$ with at least order of accuracy $\min{(\rho+1,s)}$. 
	
	\noindent By Theorem 2.3.1 in~\cite{MaWe2001}, the associated discrete Hamiltonian map has the same order of accuracy.
\end{theorem}
\begin{theorem}
	
	Suppose the Hamiltonian $H$ and its partial derivative $\frac{\partial H}{\partial p}$ are Lipschitz continuous in both variables, and $H$ is sufficiently regular for the Taylor method $\Psi_h^{(\rho+1)}$ to be well-defined. 
	
	\noindent Then $H_d^+(q_0,p_1)$ approximates $H_d^{+,E}(q_0,p_1)$ with at least order of accuracy $\min{(\rho+1,s)}$. 
	
	\noindent By Theorem 2.2 in~\cite{ScLe2017}, the associated discrete right Hamilton's map has the same order of accuracy.
\end{theorem}

Note that analogous constructions and error analysis results have been derived in~\cite{duruisseaux2020adaptive,ScShLe2017} for discrete left Hamiltonians in the Type~III case.

\subsection{Time-adaptive Hamiltonian integrators via the Poincar\'e transformation} \label{section: Poincare Hamiltonian}

Symplectic integrators form a class of geometric numerical integrators of interest since, when applied to conservative Hamiltonian systems, they yield discrete approximations of the flow that preserve the symplectic 2-form~\cite{HaLuWa2006}, which results in the preservation of many qualitative aspects of the underlying system. In particular, symplectic integrators exhibit excellent long-time near-energy preservation. However, when symplectic integrators were first used in combination with variable time-steps, the near-energy preservation was lost and the integrators performed poorly~\cite{CalSan93, GlaDunCan91}. Backward error analysis provided justification both for the excellent long-time near-energy preservation of symplectic integrators and for the poor performance experienced when using variable time-steps (see Chapter~IX of~\cite{HaLuWa2006}): it shows that symplectic integrators can be associated with a modified Hamiltonian in the form of a formal power series in terms of the time-step. The use of a variable time-step results in a different modified Hamiltonian at every iteration, which is the source of the poor energy conservation. The Poincar\'e transformation is one way to incorporate variable time-steps in geometric integrators without losing the nice conservation properties associated with these integrators.  

Given a Hamiltonian $H(q,t,p)$, consider a desired transformation of time $t \mapsto \tau$ described by the monitor function $g(q,t,p)$ via
\begin{equation}
	\frac{dt}{d\tau} = g(q,t,p).
\end{equation}
 The time $t$ shall be referred to as the physical time, while $\tau$ will be referred to as the fictive time, and we will denote derivatives with respect to $t$ and $\tau$ by dots and apostrophes, respectively. A new Hamiltonian system is constructed using the Poincar\'e transformation,
 \begin{equation}
	\bar{H}(\bar{q},\bar{p}) = g(q,\mathfrak{q},p) \left(H(q,\mathfrak{q},p) + \mathfrak{p} \right),
\end{equation}
in the extended phase space defined by $\bar{q} = \left[\begin{smallmatrix} q \\ \mathfrak{q} \end{smallmatrix} \right] \in \bar{\mathcal{Q}}$ and $\bar{p} = \left[ \begin{smallmatrix} p \\ \mathfrak{p} \end{smallmatrix} \right] $ where $\mathfrak{p}$ is the conjugate momentum for $\mathfrak{q}=t$ with
 $\mathfrak{p}(0)=-H(q(0),0,p(0))$. 
 The corresponding equations of motion in the extended phase space are then given by
 \begin{equation}
 	\bar{q} '  = \frac{\partial \bar{H}}{\partial \bar{p}},  \qquad\qquad  \bar{p}' =  -\frac{\partial \bar{H}}{\partial \bar{q}}.    
 \end{equation}
 Suppose $(\bar{Q}(\tau) , \bar{P}(\tau))$ are solutions to these extended equations of motion, and let $(q(t),p(t))$ solve Hamilton's equations for the original Hamiltonian $H$. Then
 \begin{equation}
 	\bar{H}(\bar{Q}(\tau) , \bar{P}(\tau))  = \bar{H}(\bar{Q}(0) , \bar{P}(0)) = 0.
 \end{equation}
 Therefore, the components $(Q(\tau) , P(\tau))$ in the original phase space of the augmented solutions $(\bar{Q}(\tau) , \bar{P}(\tau))$ satisfy
 \begin{equation}H(Q(\tau) , \tau , P(\tau)) = - \mathfrak{p}(\tau), \qquad  H(Q(0) , 0, P(0)) = - \mathfrak{p}(0) = H(q(0),0,p(0)) . \end{equation}
 Then,  $(Q(\tau) , P(\tau))$ and $(q(t),p(t))$ both satisfy Hamilton's equations for the original Hamiltonian~$H$ with the same initial values, so they must be the same. Note that the Hessian is given by
 \begin{equation}
 	\frac{\partial^2 \bar{H}}{\partial \bar{p}^2} = \begin{bmatrix} \frac{\partial H}{\partial p}\nabla_pg(\bar{q},p)^\top +g(\bar{q},p)\frac{\partial^2 H}{\partial p^2}+\nabla_pg(\bar{q},p)\frac{\partial H}{\partial p}^\top  \ \ \nabla_pg(\bar{q},p) \\ \nabla_pg(\bar{q},p)^\top  \qquad \qquad \qquad \qquad \qquad \qquad \qquad \qquad 0
 	\end{bmatrix},
 \end{equation}
which will be singular in many cases. The degeneracy of the Hamiltonian $\bar{H}$ implies that there is no corresponding Type~I Lagrangian formulation. This approach works seamlessly with the existing methods and theorems for Hamiltonian variational integrators, but where the system under consideration is the transformed Hamiltonian system resulting from the Poincar\'e transformation. We can use a symplectic integrator with constant time-step in fictive time $\tau$ on the Poincar\'e transformed system, which will have the effect of integrating the original system with the desired variable time-step in physical time $t$ via the relation $\frac{dt}{d\tau} = g(q,t,p)$. \\

\section{Time-adaptive Lagrangian Integrators}

The Poincar\'e transformation for time-adaptive symplectic integrators on the Hamiltonian side presented in Section~\ref{section: Poincare Hamiltonian} with autonomous monitor function $g(q,p)$ was first introduced in~\cite{Zare1975}, and extended to the case where $g$ can also depend on time based on ideas from~\cite{Ha1997}. All the literature to date on the Poincar\'e transformation have constructed the Poincar\'e transformed system by reverse-engineering: the Poincar\'e transformed Hamiltonian is chosen in such a way that the corresponding component dynamics satisfy Hamilton's equations in the original space.

\subsection{Variational Derivation of the Poincar\'e Hamiltonian} \label{section: Variational Derivation of Poincare Hamiltonian}

We now depart from the traditional reverse-engineering strategy for the Poincar\'e transformation and present a new way to think about the Poincar\'e transformed Hamiltonian by deriving it from a variational principle. This simple derivation gives additional insight into the transformation mechanism and provides natural candidates for time-adaptivity on the Lagrangian side and for more general frameworks.

As before, we work in the extended space $(q, \mathfrak{q}, p, \mathfrak{p} )$ where $\mathfrak{q} = t$ and $\mathfrak{p}$ is the corresponding conjugate momentum, and consider a time transformation $t \rightarrow \tau$ given by
\begin{equation}
	\frac{dt}{d\tau} = g(q,t,p).
\end{equation}
We define an extended action functional $\mathfrak{S} : C^2([0,T],T^*\bar{\mathcal{Q}}) \rightarrow \mathbb{R}$ by
\begin{align}
	\mathfrak{S} (\bar{q}(\cdot),\bar{p}(\cdot)) & = \bar{p}(T) \bar{q}(T) - \int_{0}^{T}{ \left[ \bar{p}(t) \dot{\bar{q}}(t) - H(q(t),t,p(t)) - \mathfrak{p}(t) \right] dt} \\ & =  \bar{p}(T) \bar{q}(T) - \int_{\tau(t=0)}^{\tau(t=T)}{ \left[ \bar{p}(\tau) \frac{d\tau}{dt} \bar{q}'(\tau) - H(q(\tau),\mathfrak{q}(\tau),p(\tau)  ) - \mathfrak{p}(\tau) \right] \frac{dt}{d\tau} d\tau } \\ & =  \bar{p}(T) \bar{q}(T) - \int_{\tau(t=0)}^{\tau(t=T)}{ \left\{ \bar{p}(\tau)  \bar{q}'(\tau) -  \frac{dt}{d\tau} \left[  H(q(\tau),\mathfrak{q}(\tau),p(\tau)  ) + \mathfrak{p}(\tau) \right]  \right\} d\tau }  ,
\end{align}
where we have performed a change of variables in the integral. Then,
\small 
\begin{equation} \label{eq: Transformed Hamiltonian Action}
	\mathfrak{S} (\bar{q}(\cdot),\bar{p}(\cdot)) =  \bar{p}(T) \bar{q}(T) - \int_{\tau(t=0)}^{\tau(t=T)}{ \left\{ \bar{p}(\tau) \bar{q}'(\tau) - g(q(\tau),\mathfrak{q}(\tau),p(\tau)) \left[ H(q(\tau),\mathfrak{q}(\tau),p(\tau)) + \mathfrak{p}(\tau) \right] \right\}  d\tau }.  
\end{equation} 
\normalsize 
Computing the variation of $\mathfrak{S}$ yields
\begin{align*} 
	\delta \mathfrak{S} & = \bar{q}(T) \delta \bar{p}(T) + \bar{p}(T) \delta \bar{q} (T)    - \int_{\tau(t=0)}^{\tau(t=T)}{ \left[ \mathfrak{q}' \delta \mathfrak{p}  + \mathfrak{p} \delta \mathfrak{q}' - \left( g \frac{\partial H}{\partial \mathfrak{q}} + \frac{\partial g}{\partial \mathfrak{q}} (H + \mathfrak{q})  \right)\delta \mathfrak{q} -  g \delta \mathfrak{p}  \right]  d\tau  }  \\  & \qquad \qquad  - \int_{\tau(t=0)}^{\tau(t=T)}{ \left[  q'  \delta p + p \delta q'  -\left( g \frac{\partial H}{\partial q} + \frac{\partial g}{\partial q} (H + \mathfrak{p})  \right) \delta q   -\left( g \frac{\partial H}{\partial p} + \frac{\partial g}{\partial p} (H + \mathfrak{p} ) \right) \delta p  \right]  d\tau  } ,
\end{align*}  
and using integration by parts and the boundary conditions $\delta \bar{q}(0) = 0$ and $\delta \bar{p} (T) = 0,$ gives
\begin{align*}
	\delta \mathfrak{S} & =  \int_{\tau(t=0)}^{\tau(t=T)}{ \left[ p' + g \frac{\partial H}{\partial q} + \frac{\partial g}{\partial q} (H + \mathfrak{p})  \right] \delta q  d\tau  } + \int_{\tau(t=0)}^{\tau(t=T)}{ \left[ g \frac{\partial H}{\partial p} + \frac{\partial g}{\partial p} (H + \mathfrak{p} )   - q' \right]  \delta p d\tau  } \\ & \qquad  \qquad + \int_{\tau(t=0)}^{\tau(t=T)}{ \left[ \mathfrak{p}' + g \frac{\partial H}{\partial \mathfrak{q}} + \frac{\partial g}{\partial \mathfrak{q}} (H + \mathfrak{p} ) \right]  \delta \mathfrak{q} d\tau  } - \int_{\tau(t=0)}^{\tau(t=T)}{ \left[ \mathfrak{q}' -g \right]  \delta \mathfrak{p} d\tau  } . 
\end{align*}  
Thus, the condition that $\mathfrak{S} (\bar{q}(\cdot),\bar{p}(\cdot)) $ is stationary with respect to the boundary conditions $\delta \bar{q}(0) = 0$ and $\delta \bar{p}(T) =0$ is equivalent to $(\bar{q}(\cdot),\bar{p}(\cdot))$ satisfying Hamilton's canonical equations corresponding to the Poincar\'e transformed Hamiltonian, 
\begin{align} 
	& 	 \mathfrak{q}' = g(q,\mathfrak{q},p), \\
	& q' = g(q,\mathfrak{q},p)  \frac{\partial H}{\partial p}(q,\mathfrak{q},p)+ \frac{\partial g}{\partial p}(q,\mathfrak{q},p)  \left[  H(q,\mathfrak{q},p) + \mathfrak{p} \right], 
	\\
	&
	p' = - g(q,\mathfrak{q},p) \frac{\partial H}{\partial q}(q,\mathfrak{q},p) - \frac{\partial g}{\partial q}(q,\mathfrak{q},p) \left[ H(q,\mathfrak{q},p)  + \mathfrak{p} \right], \\
	&	\mathfrak{p}' = - g(q,\mathfrak{q},p) \frac{\partial H}{\partial \mathfrak{q}}(q,\mathfrak{q},p) - \frac{\partial g}{\partial \mathfrak{q}}(q,\mathfrak{q},p) \left[ H(q,\mathfrak{q},p) + \mathfrak{p} \right]. 
\end{align}

An alternative way to reach the same conclusion is by interpreting equation \eqref{eq: Transformed Hamiltonian Action} as the usual Type~II action functional for the modified Hamiltonian,
\begin{equation}  g(q(\tau),\mathfrak{q}(\tau),p(\tau)) \left[ H(q(\tau),\mathfrak{q}(\tau),p(\tau)) + \mathfrak{p}(\tau) \right], \end{equation}
which coincides with the Poincar\'e transformed Hamiltonian.  \\

\subsection{Time-adaptivity from a Variational Principle on the Lagrangian side}  \label{section: Variational Time-Adaptivity Lagrangian} 

We will now derive a mechanism for time-adaptivity on the Lagrangian side by mimicking the derivation of the Poincar\'e Hamiltonian. We will work in the extended space $\bar{q} = (q,\mathfrak{q} , \lambda )^\top \in \bar{\mathcal{Q}}$ where $\mathfrak{q} = t$ and $\lambda$ is a Lagrange multiplier used to enforce the time rescaling $	\frac{dt}{d\tau} = g(t).$ Consider the action functional $\mathfrak{S} : C^2([0,T],T\bar{\mathcal{Q}} ) \rightarrow \mathbb{R}$ given by 
\begin{align}
	\mathfrak{S} (\bar{q}(\cdot),\dot{\bar{q}}(\cdot)) & =  \int_{0}^{T}{ \left[ L(q(t),\dot{q}(t),\mathfrak{q}(t)) - \lambda(t) \left( \frac{d \mathfrak{q}}{d\tau} - g(\mathfrak{q}(t)) \right) \right]dt}  \\ & =  \int_{\tau(t=0)}^{\tau(t=T)}{  \left[  \frac{dt}{d\tau} L\left(q(\tau), \frac{d\tau}{dt} q'(\tau), \mathfrak{q}(\tau) \right) - \lambda(\tau) \frac{dt}{d\tau} \left( \frac{d \mathfrak{q}}{d\tau} - g(\mathfrak{q}(\tau)) \right) \right]   d\tau }  
	\\ & =  \int_{\tau(t=0)}^{\tau(t=T)}{  \left[  \mathfrak{q}'(\tau) L\left(q(\tau), \frac{d\tau}{dt} q'(\tau), \mathfrak{q}(\tau) \right) - \lambda(\tau)\mathfrak{q}'(\tau) \left[ \mathfrak{q}'(\tau)- g(\mathfrak{q}(\tau)) \right]   \right]  d\tau } , 
\end{align}
where, as before, we have performed a change of variables in the integral. This is the usual Type~I action functional for the extended autonomous Lagrangian,
\begin{equation} \bar{L} (\bar{q}(\tau),\bar{q}'(\tau) ) =  \mathfrak{q}'(\tau) L\left(q(\tau), \frac{d\tau}{dt} q'(\tau), \mathfrak{q}(\tau) \right) - \lambda(\tau)\mathfrak{q}'(\tau) \left[ \mathfrak{q}'(\tau)- g(\mathfrak{q}(\tau)) \right] . \end{equation}

\begin{theorem}
	If $\left(\bar{q}(\tau), \bar{q}'(\tau) \right)$ satisfies the Euler--Lagrange equations corresponding to the extended Lagrangian $\bar{L}$, then its components satisfy $\frac{dt}{d\tau} = g(t)$ and the original Euler--Lagrange equations 
	\begin{equation}
		\frac{d}{dt}  \frac{\partial L }{\partial \dot{q}}  \left(q, \dot{q}, t\right)  =   \frac{\partial L}{\partial q} \left(q, \dot{q} , t \right)   .
	\end{equation}
	\proof{
		Substituting the expression for $\bar{L}$ into the Euler--Lagrange equations, $ \frac{d}{d\tau} \frac{\partial \bar{L}}{\partial \lambda'}  = \frac{\partial \bar{L}}{\partial \lambda}, $ and $   \frac{d}{d\tau} \frac{\partial \bar{L}}{\partial q'}  = \frac{\partial \bar{L}}{\partial q}$, 		gives  \[  \mathfrak{q}' \left[ \mathfrak{q}' - g(\mathfrak{q}) \right] = 0 ,\] and
		\[ \frac{d\mathfrak{q}}{d\tau} \frac{d}{d\mathfrak{q}} \left[ \mathfrak{q}' \frac{\partial L \left(q, \frac{d\tau}{d\mathfrak{q}} q'  , \mathfrak{q} \right) }{\partial q'} \right] =  \mathfrak{q}'  \frac{\partial L \left(q, \frac{d\tau}{d\mathfrak{q}} q', \mathfrak{q} \right)}{\partial q} . \] 
		Now, $\mathfrak{q}' = g(\mathfrak{q}) >0 $ so $\mathfrak{q}' = g(\mathfrak{q})$, and the chain rule gives
		\[ \frac{d}{d\mathfrak{q}}  \frac{\partial L }{\partial \dot{q}} \left(q, \frac{d\tau}{d\mathfrak{q}} q' , \mathfrak{q} \right)  = \frac{\partial L }{\partial q} \left(q, \frac{d\tau}{d\mathfrak{q}} q' , \mathfrak{q} \right)  .  \]  
		Using the equation $\dot{q}  =\frac{d\tau}{d\mathfrak{q}} q'  $ and replacing $\mathfrak{q}$ by $t$ recovers the original Euler--Lagrange equations. \qed \\}
\end{theorem}

We now introduce a discrete variational formulation of these continuous Lagrangian mechanics. Suppose we are given a partition $0 = \tau_0 < \tau_1 < \ldots < \tau_N = \mathcal{T}$ of the interval $[0,\mathcal{T}]$, and a discrete curve in $\mathcal{Q} \times \mathbb{R} \times \mathbb{R} $ denoted by $\{  (q_k, \mathfrak{q}_k, \lambda _k) \}_{k=0}^{N}$ such that $q_k \approx q(\tau_k)$,   $\mathfrak{q}_k \approx \mathfrak{q}(\tau_k)$, and $\lambda_k \approx \lambda(\tau_k)$. Consider the discrete action functional,
\begin{equation} 
	\bar{\mathfrak{S}}_d \left(\{  (q_k, \mathfrak{q}_k, \lambda _k) \}_{k=0}^{N} \right)  =  \sum_{k=0}^{N-1}{  \left[  L_d(q_k, \mathfrak{q}_k,  q_{k+1} , \mathfrak{q}_{k+1}) -  \lambda_{k} \frac{ \mathfrak{q}_{k+1} - \mathfrak{q}_k}{\tau_{k+1} - \tau_{k}}  + \lambda_k g(\mathfrak{q}_k)   \right]  \frac{ \mathfrak{q}_{k+1} - \mathfrak{q}_k}{\tau_{k+1} - \tau_{k}}    },
\end{equation}
where $L_d(q_k, \mathfrak{q}_k,  q_{k+1} , \mathfrak{q}_{k+1}) $ is obtained by approximating the exact discrete Lagrangian, which is related to Jacobi's solution of the Hamilton--Jacobi equation and is the generating function for the exact time-$h$ flow map. It is given by the extremum of the action integral from $\tau_k$ to $\tau_{k+1}$ over twice continuously differentiable curves $(q,\mathfrak{q}) \in \mathcal{Q}\times \mathbb{R}$ satisfying the boundary conditions $ (q(\tau_k), \mathfrak{q}(\tau_k)) = (q_k, \mathfrak{q}_k)  ,$ and $ (q(\tau_{k+1}), \mathfrak{q}(\tau_{k+1})) = (q_{k+1}, \mathfrak{q}_{k+1})$:
\begin{equation}
	L_d(q_k, \mathfrak{q}_k,  q_{k+1} , \mathfrak{q}_{k+1})  \approx  \ext_{ \substack{ (q,\mathfrak{q} ) \in C^2 ([\tau_{k} , \tau_{k+1} ] , \mathcal{Q} \times \mathbb{R} ) \\ (q, \mathfrak{q})(\tau_k) = (q_k, \mathfrak{q}_k)  , \text{  } (q, \mathfrak{q})(\tau_{k+1}) = (q_{k+1}, \mathfrak{q}_{k+1}) }  }   \text{  } \int_{\tau_k}^{\tau_{k+1}}{L\left(q, \frac{q'}{g(\mathfrak{q})}  , \mathfrak{q} \right) d\tau }.
\end{equation}
In practice, we can obtain an approximation by replacing the integral with a quadrature rule, and extremizing over a finite-dimensional function space instead of $C^2 ([\tau_{k} , \tau_{k+1} ] , \mathcal{Q} \times \mathbb{R} )$.
This discrete functional $\bar{\mathfrak{S}}_d$ is a discrete analogue of the action functional $\bar{\mathfrak{S}} : C^2([0,T], \mathcal{Q} \times \mathbb{R} \times \mathbb{R}) \rightarrow \mathbb{R}$ given by
\begin{align} 
	\bar{\mathfrak{S}} (q(\cdot), \mathfrak{q}(\cdot), \lambda(\cdot)) & =  \int_{0}^{\mathcal{T}}{ \bar{L} \left(q(\tau), \mathfrak{q}(\tau),\lambda(\tau),q'(\tau), \mathfrak{q}'(\tau) , \lambda'(\tau) \right)   d\tau}  \\ & =    \int_{0}^{\mathcal{T}}{  \left[ L\left(q, \frac{q'}{g(\mathfrak{q})}  , \mathfrak{q} \right) -\lambda \mathfrak{q}' + \lambda g(\mathfrak{q}) \right] \mathfrak{q}' d\tau}  .
\end{align}
\hfill  \\

We can derive the following result which relates a discrete Type~I variational principle to a set of discrete Euler--Lagrange equations:

\begin{theorem}  \label {Theorem: discrete EL equations 2}The Type~I discrete Hamilton's variational principle,
	\begin{equation}
		\delta \bar{\mathfrak{S}}_d \left(\{  (q_k, \mathfrak{q}_k, \lambda _k) \}_{k=0}^{N} \right)  = 0,
	\end{equation}
	is equivalent to the discrete extended Euler--Lagrange equations,
	\begin{equation}   \mathfrak{q}_{k+1} =  \mathfrak{q}_k + (\tau_{k+1} - \tau_{k} )g(\mathfrak{q}_k)  ,\end{equation}  
	\begin{equation}   \frac{ \mathfrak{q}_{k+1} - \mathfrak{q}_k}{\tau_{k+1} - \tau_{k}}   D_1 L_d(q_k, \mathfrak{q}_k,  q_{k+1} , \mathfrak{q}_{k+1}) + \frac{ \mathfrak{q}_{k} - \mathfrak{q}_{k-1}}{\tau_{k} - \tau_{k-1}}   D_3 L_d(q_{k-1}, \mathfrak{q}_{k-1},  q_{k} , \mathfrak{q}_{k})     = 0, 	\end{equation}  
	\small
\begin{equation}
	\begin{aligned}
		&    \left[    D_2 L_{d_k} + \lambda_{k} \frac{ 1}{\tau_{k+1} - \tau_{k}}  + \lambda_k \nabla  g(\mathfrak{q}_k)  \right] \frac{ \mathfrak{q}_{k+1} - \mathfrak{q}_k}{\tau_{k+1} - \tau_{k}}      -  \frac{ 1}{\tau_{k+1} - \tau_{k}} \left[    L_{d_k} -  \lambda_{k} \frac{ \mathfrak{q}_{k+1} - \mathfrak{q}_k}{\tau_{k+1} - \tau_{k}}  + \lambda_k g(\mathfrak{q}_k)      \right]  \\ & \qquad    + \left[    D_4 L_{d_{k-1}} - \lambda_{k-1} \frac{ 1}{\tau_{k} - \tau_{k-1}}   \right] \frac{ \mathfrak{q}_{k} - \mathfrak{q}_{k-1}}{\tau_{k} - \tau_{k-1}}   +  \frac{ 1}{\tau_{k} - \tau_{k-1}} \left[     L_{d_{k-1}} -  \lambda_{k-1} \frac{ \mathfrak{q}_{k} - \mathfrak{q}_{k-1}}{\tau_{k} - \tau_{k-1}}  + \lambda_{k-1} g(\mathfrak{q}_{k-1})      \right]     =0,
	\end{aligned} 
\end{equation}
\normalsize 	where $L_{d_k} $ denotes $ L_d(q_k, \mathfrak{q}_k,  q_{k+1} , \mathfrak{q}_{k+1})$.
	\proof{See Appendix \ref{Appendix: Discrete Proof 1}.		\qed \\}
\end{theorem}  

Defining the discrete momenta via the discrete Legendre transformations,
\begin{equation} p_k = - D_1 L_d(q_k,\mathfrak{q}_k,q_{k+1},\mathfrak{q}_{k+1}),  \qquad \qquad  \mathfrak{p}_k = - D_2 L_d(q_k,\mathfrak{q}_k,q_{k+1},\mathfrak{q}_{k+1}), \end{equation}
and using a constant time-step $h$ in $\tau$, the discrete Euler--Lagrange equations can be rewritten as
\begin{align} p_k  & = - D_1 L_d(q_k,\mathfrak{q}_k,q_{k+1},\mathfrak{q}_{k+1}) , \\  \mathfrak{p}_k & = - D_2 L_d(q_k,\mathfrak{q}_k,q_{k+1},\mathfrak{q}_{k+1}) , \\    \mathfrak{q}_{k+1} & =  \mathfrak{q}_k + hg(\mathfrak{q}_k) ,  \\
	p_{k+1}    & = \frac{g(\mathfrak{q}_k)}{g(\mathfrak{q}_{k+1})}  D_3 L_d(q_k, \mathfrak{q}_k,  q_{k+1} , \mathfrak{q}_{k+1})   ,  \\ \mathfrak{p}_{k+1} & =    \frac{ L_{d_k}  -  L_{d_{k+1}}}{ hg( \mathfrak{q}_{k+1} )}    + \frac{ \lambda_{k+1}}{h}  + \lambda_{k+1} \nabla  g(\mathfrak{q}_{k+1}) + \frac{g(\mathfrak{q}_k)}{g(\mathfrak{q}_{k+1})}    \left[    D_4 L_{d_k} - \frac{ \lambda_{k} }{h}   \right]  .  \end{align}  \\

\subsection{A Second Time-Adaptive Framework obtained by Reverse-Engineering}  \label{section: Guess Time-Adaptivity Lagrangian}

As mentioned earlier, all the literature to date on the Poincar\'e transformation have constructed the Poincar\'e transformed system by reverse-engineering. The Poincar\'e transformed Hamiltonian is chosen in such a way that the corresponding component dynamics satisfy the Hamilton's equations in the original space. We will follow a similar strategy to derive a second framework for time-adaptivity from the Lagrangian perspective. 

Given a time-dependent Lagrangian $L(q(t),\dot{q}(t),t)$ consider a transformation of time $t \rightarrow \tau$,
\begin{equation}
	\frac{dt}{d\tau} = g(t),
\end{equation}
described by the monitor function $g(t)$. The time $t$ shall be referred to as the physical time, while $\tau$ will be referred to as the fictive time, and we will denote derivatives with respect to $t$ and $\tau$ by dots and apostrophes, respectively. We define the autonomous Lagrangian,
\begin{equation} \label{eq: Transformed L} \bar{L} (\bar{q}(\tau),\bar{q}'(\tau) ) = \mathfrak{q}' L\left(q, \frac{q'}{g(\mathfrak{q})}  , \mathfrak{q} \right)  - \lambda \left( \mathfrak{q}' - g(\mathfrak{q}) \right) , \end{equation}
in the extended space with $\bar{q} = \left( q , \mathfrak{q}, \lambda  \right)^\top$ where $\mathfrak{q} = t$, and where $\lambda$ is a multiplier used to impose the constraint that  the time evolution is guided by the monitor function $g(t)$. Note that in contrast to the earlier framework, the Lagrange multiplier term lacks an extra multiplicative factor of $\mathfrak{q}'$.

\begin{theorem}
	If $\left(\bar{q}(\tau), \bar{q}'(\tau) \right)$ satisfies the Euler--Lagrange equations corresponding to the extended Lagrangian $\bar{L}$, then its components satisfy $\frac{dt}{d\tau} = g(t)$ and the original Euler--Lagrange equations
	\begin{equation}
		\frac{d}{dt}  \frac{\partial L }{\partial \dot{q}}  \left(q, \dot{q}, t\right)  =   \frac{\partial L}{\partial q} \left(q, \dot{q} , t \right)   .
	\end{equation}
	\begin{proof}
		Substituting the expression for $\bar{L}$ in the Euler--Lagrange equations $\frac{d}{d\tau} \frac{\partial \bar{L}}{\partial \lambda'}  = \frac{\partial \bar{L}}{\partial \lambda} $ and $   \frac{d}{d\tau} \frac{\partial \bar{L}}{\partial q'}  = \frac{\partial \bar{L}}{\partial q} $
		 gives  \[\mathfrak{q}' = g(\mathfrak{q}),\] and
		\[ \frac{d\mathfrak{q}}{d\tau} \frac{d}{d\mathfrak{q}} \left[ \frac{d\mathfrak{q}}{d\tau}  \frac{\partial L \left(q, \frac{q'}{g(\mathfrak{q})}  , \mathfrak{q} \right) }{\partial q'} \right] =  \ \frac{d\mathfrak{q}}{d\tau}  \frac{\partial L \left(q, \frac{q'}{g(\mathfrak{q})}  , \mathfrak{q} \right)}{\partial q} . \] 
			We can divide by $\frac{d\mathfrak{q}}{d\tau}$ and use the chain rule to get $\mathfrak{q}' = g(\mathfrak{q})$ and
		\[ \frac{d}{d\mathfrak{q}}  \frac{\partial L }{\partial \dot{q}} \left(q, \frac{d\tau}{d\mathfrak{q}} q' , \mathfrak{q} \right)  = \frac{\partial L }{\partial q} \left(q, \frac{d\tau}{d\mathfrak{q}} q' , \mathfrak{q} \right)  .  \]  
		Using the equations $\dot{q}  =\frac{d\tau}{d\mathfrak{q}} q'  $ and $ \mathfrak{q}' = g(\mathfrak{q}) $, and replacing $\mathfrak{q}$ by $t$ recovers the desired equations.
	\end{proof}
\end{theorem}

We now introduce a discrete variational formulation of these continuous Lagrangian mechanics. Suppose we are given a partition $0 = \tau_0 < \tau_1 < \ldots < \tau_N = \mathcal{T}$ of the interval $[0,\mathcal{T}]$, and a discrete curve in $\mathcal{Q} \times \mathbb{R} \times \mathbb{R} $ denoted by $\{  (q_k, \mathfrak{q}_k, \lambda _k) \}_{k=0}^{N}$ such that $q_k \approx q(\tau_k)$,   $\mathfrak{q}_k \approx \mathfrak{q}(\tau_k)$, and $\lambda_k \approx \lambda(\tau_k)$. Consider the discrete action functional,
\begin{equation} 
	\bar{\mathfrak{S}}_d \left(\{  (q_k, \mathfrak{q}_k, \lambda _k) \}_{k=0}^{N} \right)  =  \sum_{k=0}^{N-1}{\left\{ \frac{ \mathfrak{q}_{k+1} - \mathfrak{q}_k}{\tau_{k+1} - \tau_{k}}  \left[  L_d(q_k, \mathfrak{q}_k,  q_{k+1} , \mathfrak{q}_{k+1}) -  \lambda_{k} \right] + \lambda_k g(\mathfrak{q}_k)     \right\}},
\end{equation}
where,
\begin{equation}
	L_d(q_k, \mathfrak{q}_k,  q_{k+1} , \mathfrak{q}_{k+1})  \approx  \ext_{ \substack{ (q,\mathfrak{q} ) \in C^2 ([\tau_{k} , \tau_{k+1} ] , \mathcal{Q} \times \mathbb{R} ) \\ (q, \mathfrak{q})(\tau_k) = (q_k, \mathfrak{q}_k) , \text{  }  (q, \mathfrak{q})(\tau_{k+1}) = (q_{k+1}, \mathfrak{q}_{k+1}) }  }   \text{  } \int_{\tau_k}^{\tau_{k+1}}{L\left(q, \frac{q'}{g(\mathfrak{q})}  , \mathfrak{q} \right) d\tau }.
\end{equation}
This discrete functional $\bar{\mathfrak{S}}_d$ is a discrete analogue of the action functional $\bar{\mathfrak{S}} : C^2([0,T], \mathcal{Q} \times \mathbb{R} \times \mathbb{R}) \rightarrow \mathbb{R}$ given by
\begin{align} 
	\bar{\mathfrak{S}} (q(\cdot), \mathfrak{q}(\cdot), \lambda(\cdot)) & =  \int_{0}^{\mathcal{T}}{ \bar{L} \left(q(\tau), \mathfrak{q}(\tau),\lambda(\tau),q'(\tau), \mathfrak{q}'(\tau) , \lambda'(\tau) \right)   d\tau}  \\ & =    \int_{0}^{\mathcal{T}}{ \left\{    \mathfrak{q}' \left[ L\left(q, \frac{q'}{g(\mathfrak{q})}  , \mathfrak{q} \right) -\lambda  \right] + \lambda  g(\mathfrak{q})  \right\}  d\tau}  .
\end{align}

We can derive the following result which relates a discrete Type~I variational principle to a set of discrete Euler--Lagrange equations:

\begin{theorem} \label{Theorem: discrete EL equations 1} The Type~I discrete Hamilton's variational principle,
	\begin{equation}
		\delta \bar{\mathfrak{S}}_d \left(\{  (q_k, \mathfrak{q}_k, \lambda _k) \}_{k=0}^{N} \right)  = 0,
	\end{equation}
	is equivalent to the discrete extended Euler--Lagrange equations,
			\begin{equation}  \mathfrak{q}_{k+1} =  \mathfrak{q}_k + (\tau_{k+1} - \tau_{k} )g(\mathfrak{q}_k),  \end{equation}  
		\begin{equation}  \frac{ \mathfrak{q}_{k+1} - \mathfrak{q}_k}{\tau_{k+1} - \tau_{k}}  D_1 L_d(q_k, \mathfrak{q}_k,  q_{k+1} , \mathfrak{q}_{k+1}) +  \frac{ \mathfrak{q}_{k} - \mathfrak{q}_{k-1}}{\tau_{k} - \tau_{k-1}}  D_3 L_d(q_{k-1}, \mathfrak{q}_{k-1},  q_{k} , \mathfrak{q}_{k})  = 0, \end{equation} 
		\small 
	\begin{align}
		&  \frac{ \mathfrak{q}_{k+1} - \mathfrak{q}_k}{\tau_{k+1} - \tau_{k}}  D_2 L_{d_k} -  \frac{ L_{d_k}}{\tau_{k+1} - \tau_{k}}+ \frac{ \mathfrak{q}_{k} - \mathfrak{q}_{k-1}}{\tau_{k} - \tau_{k-1}}  D_4 L_{d_{k-1}} + \frac{  L_{d_{k-1} }}{\tau_{k} - \tau_{k-1}}  =  \frac{\lambda_{k-1}}{\tau_{k} - \tau_{k-1}}  - \frac{\lambda_k}{\tau_{k+1} - \tau_k}  - \lambda_k  \nabla  g(\mathfrak{q}_k)    ,
	\end{align}
	where $L_{d_k} $ denotes $ L_d(q_k, \mathfrak{q}_k,  q_{k+1} , \mathfrak{q}_{k+1})$.
	\proof{See Appendix \ref{Appendix: Discrete Proof 2}.		\qed}
\end{theorem}  

\normalsize

Defining the discrete momenta via the discrete Legendre transformations,
\begin{equation} p_k = - D_1 L_d(q_k,\mathfrak{q}_k,q_{k+1},\mathfrak{q}_{k+1}),  \qquad \qquad  \mathfrak{p}_k = - D_2 L_d(q_k,\mathfrak{q}_k,q_{k+1},\mathfrak{q}_{k+1}), \end{equation}
and using a constant time-step $h$ in $\tau$, the discrete Euler--Lagrange equations can be rewritten as
\begin{align} p_k  & = - D_1 L_d(q_k,\mathfrak{q}_k,q_{k+1},\mathfrak{q}_{k+1}) , \\  \mathfrak{p}_k & = - D_2 L_d(q_k,\mathfrak{q}_k,q_{k+1},\mathfrak{q}_{k+1}) , \\    \mathfrak{q}_{k+1} & =  \mathfrak{q}_k + hg(\mathfrak{q}_k)  ,\\
	p_{k+1}    & = \frac{g(\mathfrak{q}_k)}{g(\mathfrak{q}_{k+1})}  D_3 L_d(q_k, \mathfrak{q}_k,  q_{k+1} , \mathfrak{q}_{k+1}) ,  \\ \mathfrak{p}_{k+1} & = \frac{L_{d_k}  -   L_{d_{k+1}} + \lambda_{k+1} -\lambda_k  + h\lambda_{k+1}  \nabla g(\mathfrak{q}_{k+1})  + h g(\mathfrak{q}_k) D_4L_{d_k}  }{hg(\mathfrak{q}_{k+1})} .   \end{align}  \\

\subsection{Remarks on the Framework for Time-Adaptivity} Time-adaptivity comes more naturally on the Hamiltonian side through the Poincar\'e transformation. Indeed, in the Hamiltonian case, the time-rescaling equation $ \mathfrak{q}' = g(q,\mathfrak{q},p)$ emerged naturally through the change of time variable inside the extended action functional. By contrast, in the Lagrangian case, we need to impose the time-rescaling equation as a constraint via a multiplier, which we then consider as an extra position coordinate. This strategy can be thought of as being part of the more general framework for constrained variational integrators (see~\cite{MaWe2001,Duruisseaux2022Constrained}). 

The Poincar\'e transformation on the Hamiltonian side was presented in 
\cite{Zare1975,Ha1997,duruisseaux2020adaptive} for the general case where the monitor function can depend on position, time and momentum, $g= g(q,t,p)$. For the accelerated optimization application which was our main motivation to develop a time-adaptive framework for geometric integrators, the monitor function only depends on time, $g = g(t)$. For the sake of simplicity and clarity, we have decided to only present the theory for time-adaptive Lagrangian integrators for monitor functions of the form $g = g(t)$ in this paper. Note however that this time-adaptivity framework on the Lagrangian side can be extended to the case where the monitor function also depends on position, $g = g(q,t)$. The action integral remains the same with the exception that $g$ is now a function of $(q,\mathfrak{q})$. Unlike the case where $g=g(t)$, the corresponding Euler--Lagrange equation $	\frac{d}{d\tau}  \frac{\partial \bar{L} }{\partial q'}  =   \frac{\partial L}{\partial q}$ yields an extra term $\lambda(t)  \frac{\partial g}{\partial q} (q,t)$ in the original phase-space,
\begin{equation} \frac{d}{dt} \frac{\partial L}{\partial \dot{q}}  \left(q, \dot{q}, t \right)  - \frac{\partial L}{\partial q} \left(q, \dot{q}, t\right) =    \lambda(t)  \frac{\partial g}{\partial q} (q,t). \end{equation}
The discrete Euler--Lagrange equations become more complicated and involve terms with partial derivatives  $D_1 g(q_k,\mathfrak{q}_k)$ of $g$ with respect to $q$. Furthermore, when $g = g(q,t)$, the discrete Euler--Lagrange equations involve $\lambda_k $ but the time-evolution of the Lagrange multiplier $\lambda$ is not well-defined, so the discrete Hamiltonian map corresponding to the discrete Lagrangian $L_d$ is not well-defined, as explained in~\cite[page 440]{MaWe2001}. Although there are ways to circumvent this problem in practice, this adds some difficulty and makes the time-adaptive Lagrangian approach with $g = g(q,t)$ less natural and desirable than the corresponding Poincar\'e transformation on the Hamiltonian side. It might also be tempting to generalize further and consider the case where $g = g(q,\dot{q},t)$. However, in this case, the time-rescaling equation $\frac{dt}{d\tau } = g(q,\dot{q},t)$ becomes implicit and it becomes less clear how to generalize the variational derivation presented in this paper. There are examples where time-adaptivity with these more general monitor functions proved advantageous (see for instance Kepler's problem in~\cite{duruisseaux2020adaptive}). This motivates further effort towards developing a better framework for time-adaptivity on the Lagrangian side with more general monitor functions.

It might be more natural to consider these time-rescaled Lagrangian and Hamiltonian dynamics as Dirac mechanics~\cite{LeOh2008,YoMa2006a,YoMa2006b} on the Pontryagin bundle $(q,v,p) \in T\mathcal{Q} \oplus T^* \mathcal{Q}$. Dirac dynamics are described by the Hamilton-Pontryagin variational principle where the momentum $p$ acts as a Lagrange multiplier to impose the kinematic equation $\dot{q} = v$,
\begin{equation}
	\delta \int_{0}^{T}{\left[  L(q,v,t) + p (\dot{q} - v) \right] dt} = 0.
\end{equation}
This provides a variational description of both Lagrangian and Hamiltonian mechanics, yields the
implicit Euler–Lagrange equations
\begin{equation}
	\dot{q} = v, \qquad \dot{p} = \frac{\partial L}{\partial q}, \qquad  p = \frac{\partial L}{\partial v},
\end{equation}
and suggests the introduction of a more general quantity, the generalized energy \begin{equation} E(q,v,p,t) = pv - L(q,v,t),   \end{equation}
as an alternative to the Hamiltonian.

\section{Application to Accelerated Optimization on Vector Spaces} \label{section: Accelerated Optimization Vector Spaces}

\subsection{A Variational Framework for Accelerated Optimization}

A variational framework was introduced in~\cite{WiWiJo16} for accelerated optimization on normed vector spaces. The $p$-Bregman Lagrangians and Hamiltonians are defined to be
\begin{equation}
	\mathcal{L}_{p}(x,v,t) = \frac{t^{ p +1}}{2p} \langle v , v\rangle  - Cpt^{2p-1} f(x),  \end{equation} 
\begin{equation} 
	\mathcal{H}_{p}(x,r,t)= \frac{p}{2t^{p +1}} \langle r , r\rangle + Cpt^{2p-1} f(x), 
\end{equation}
which are scalar-valued functions of position $x\in \mathcal{X}$, velocity $v\in \mathbb{R}^d$ or  momentum $r\in \mathbb{R}^d$, and time $t$. In~\cite{WiWiJo16}, it was shown that solutions to the  $p$-Bregman Euler--Lagrange equations converge to a minimizer of the convex objective function $f$ at a convergence rate of $\mathcal{O}(1/t^p)$. Furthermore, this family of Bregman dynamics is closed under time dilation: time-rescaling a solution to the $p$-Bregman Euler--Lagrange equations via $\tau(t) =t^{\mathring{p}/p}$ yields a solution to the $\mathring{p}$-Bregman Euler--Lagrange equations. Thus, the entire subfamily of Bregman trajectories indexed by the parameter $p$ can be obtained by speeding up or slowing down along the Bregman curve corresponding to any value of $p$. %The associated $p$-Bregman Euler--Lagrange equations are given by
%\begin{equation}  \ddot{x} + \frac{ p +1}{t} \dot{x} + Cp^2t^{p-2} \nabla f(x) = 0, \end{equation} 
In~\cite{duruisseaux2020adaptive}, the time-rescaling property of the Bregman dynamics was exploited together with a carefully chosen Poincar\'e transformation to transform the $p$-Bregman Hamiltonian into an autonomous version of the $\mathring{p}$-Bregman Hamiltonian in extended phase-space, where $\mathring{p} < p$. This strategy allowed us to achieve the faster rate of convergence associated with the higher-order $p$-Bregman dynamics, but with the computational efficiency of integrating the lower-order $\mathring{p}$-Bregman dynamics. Explicitly, using the time rescaling $\tau(t) = t^{\mathring{p}/p}$ within the Poincar\'e transformation framework yields the adaptive approach $p\rightarrow \mathring{p}$-Bregman Hamiltonian,
\begin{equation} 
	\begin{aligned}
		\bar{H}_{p\rightarrow \mathring{p}}(\bar{q},\bar{r}) & =   \frac{p^2}{2\mathring{p}\mathfrak{q}^{p+\mathring{p}/p} }  \langle r , r \rangle  + \frac{Cp^2}{\mathring{p}} \mathfrak{q}^{2p-\mathring{p}/p} f(q) + \frac{p}{\mathring{p}} \mathfrak{r} \mathfrak{q}^{1-\mathring{p}/p} ,
	\end{aligned}
\end{equation}
and when $\mathring{p} = p$, the direct approach $p$-Bregman Hamiltonian,
\begin{equation}
	\bar{H}_{p}(\bar{q},\bar{r}) = \frac{p}{2\mathfrak{q}^{p+1}} \langle r , r \rangle + Cp\mathfrak{q}^{2p-1} f(q) + \mathfrak{r}.
\end{equation}

In~\cite{duruisseaux2020adaptive}, a careful computational study was performed on how time-adaptivity and symplecticity of the numerical scheme improve the performance of the resulting optimization algorithm. In particular, it was observed that time-adaptive Hamiltonian variational discretizations, which are automatically symplectic, with adaptive time-steps informed by the time-rescaling of the family of $p$-Bregman Hamiltonians yielded the most robust and computationally efficient optimization algorithms, outperforming fixed-timestep symplectic discretizations, adaptive-timestep non-symplectic discretizations, and Nesterov's accelerated gradient algorithm which is neither time-adaptive nor symplectic.  \\

\subsection{Numerical Methods} \label{section: Numerical Methods}

\subsubsection{A Lagrangian Taylor Variational Integrator (LTVI)} \label{section: LTVI Algorithm}

We will now construct a time-adaptive Lagrangian Taylor variational integrator (LTVI) for the $p$-Bregman Lagrangian,
\begin{equation}    \bar{L}_p\left(q, q' , \mathfrak{q} \right) = \frac{\mathfrak{q}^{p +1}}{2p} \langle q' , q' \rangle  - Cp\mathfrak{q}^{2p-1} f(q) ,
\end{equation}
using the strategy outlined in Section~\ref{section: Variational Integrators} together with the discrete Euler--Lagrange equations derived in Sections~\ref{section: Variational Time-Adaptivity Lagrangian} and \ref{section: Guess Time-Adaptivity Lagrangian}. 

Looking at the form of the continuous $p$-Bregman Euler--Lagrange equations,
\begin{equation}
	\ddot{q} + \frac{p+1}{\mathfrak{q}} \dot{q} + Cp^2 \mathfrak{q}^{p-2} \nabla f(q) = 0,
\end{equation}
 we can note that $\nabla f$ appears in the expression for $\ddot{q}$. Now, the construction of a LTVI as presented in Section~\ref{section: Variational Integrators} requires $\rho$-order and $(\rho +1)$-order Taylor approximations of $q$. This means that if we take $\rho \geq 1$, then $\nabla f$ and higher-order derivatives of $f$ will appear in the resulting discrete Lagrangian~$L_d$, and as a consequence, the discrete Euler--Lagrange equations,  
 	\begin{equation}
 	p_0=-D_1 L_d(q_0, q_1),\qquad p_1=D_2 L_d(q_0, q_1),
 \end{equation}
will yield an integrator which is not gradient-based. Keeping in mind the machine learning applications where data sets are very large, we will restrict ourselves to explicit first-order optimization algorithms, and therefore the highest value of $\rho$ that we can choose to obtain a gradient-based algorithm is $\rho = 0$. 

With $\rho = 0$, the choice of quadrature rule does not matter, so we can take the rectangular quadrature rule about the initial point  ($c_1 = 0$ and $b_1 =1$). We first approximate $\dot{\bar{q}}(0)=\bar{v}_0$ by the solution $\tilde{\bar{v}}_0$ of the problem 
$ \bar{q}_1 = \pi_{Q} \circ \Psi_h^{(1)}(\bar{q}_0,\tilde{\bar{v}}_0) = \bar{q}_0 + h \tilde{\bar{v}}_0$, that is $\tilde{\bar{v}}_0 = \frac{\bar{q}_1 - \bar{q}_0}{h}$. Then, applying the quadrature rule gives the associated discrete Lagrangian,
\begin{align}
	L_d(\bar{q}_0,\bar{q}_1)  & = h L_{p} \left(q_0 , \frac{\tilde{v}_0}{g(\mathfrak{q}_0)} , \mathfrak{q}_0 \right) =  \frac{\mathfrak{q}_0^{p +1}}{2p (g(\mathfrak{q}_0))^2} h \langle \tilde{v}_0 , \tilde{v}_0 \rangle  - Chp\mathfrak{q}_0^{2p-1} f(q_0).
\end{align} 
The variational integrator is then defined by the discrete extended Euler--Lagrange equations derived in Sections~\ref{section: Variational Time-Adaptivity Lagrangian} and \ref{section: Guess Time-Adaptivity Lagrangian}. In practice, we are not interested in the evolution of the conjugate momentum~$\mathfrak{r}$, and since it will not appear in the updates for the other variables, the discrete equations of motion from Sections~\ref{section: Variational Time-Adaptivity Lagrangian} and \ref{section: Guess Time-Adaptivity Lagrangian} both reduce to the same updates,
\begin{align} r_k  & = - D_1 L_d(q_k,\mathfrak{q}_k,q_{k+1},\mathfrak{q}_{k+1}) ,  \\   
	r_{k+1}    & = \frac{g(\mathfrak{q}_k)}{g(\mathfrak{q}_{k+1})}  D_3 L_d(q_k, \mathfrak{q}_k,  q_{k+1} , \mathfrak{q}_{k+1}) ,   \\  \mathfrak{q}_{k+1} & =  \mathfrak{q}_k + hg(\mathfrak{q}_k).   \end{align}
Now, for the adaptive approach, substituting $g(\mathfrak{q}) = \frac{p}{\mathring{p}} \mathfrak{q}^{1-\frac{\mathring{p}}{p}}  $  and 
\begin{equation}
	L_d(q_k , \mathfrak{q}_k,q_{k+1} , \mathfrak{q}_{k+1}) = \frac{\mathring{p} ^2}{2hp^3 } \mathfrak{q}_k^{p -1 + 2 \mathring{p} / p} \langle q_{k+1} - q_k , q_{k+1}-q_k \rangle  - Chp\mathfrak{q}_0^{2p-1} f(q_k) ,
\end{equation} 
yields the adaptive LTVI algorithm,
\begin{align}   \mathfrak{q}_{k+1}  & = \mathfrak{q}_k + h \frac{p}{\mathring{p}} \mathfrak{q}_k^{1-\mathring{p}/p} , \\ q_{k+1} & = q_k +  \frac{ hp^3  }{ \mathring{p} ^2\mathfrak{q}_k^{p -1 + 2 \mathring{p} / p}} r_k - \frac{ Ch^2p^4  }{ \mathring{p} ^2} \mathfrak{q}_k^{p - 2 \mathring{p} / p } \nabla f(q_k) ,  \\  r_{k+1}     & =   \frac{\mathring{p} ^2\mathfrak{q}_k^{p  +  \mathring{p} / p}}{hp^3 \mathfrak{q}_{k+1}^{1-\mathring{p}/p} } (q_{k+1}-q_k)      .
\end{align} 

\noindent In the direct approach, $\mathring{p} =p$ so $g(\mathfrak{q}) = 1 $ and we obtain the direct LTVI algorithm,
\begin{align}   \mathfrak{q}_{k+1}  & = \mathfrak{q}_k + h , \\ q_{k+1} & = q_k +  \frac{ hp  }{ \mathfrak{q}_k^{p +1}} r_k - Ch^2p^2   \mathfrak{q}_k^{p - 2  } \nabla f(q_k) ,  \\  r_{k+1}     & =   \frac{\mathfrak{q}_k^{p  + 1}}{hp } (q_{k+1}-q_k)      . 
\end{align}  

\subsubsection{A Hamiltonian Taylor Variational Integrator (HTVI)} 

In~\cite{duruisseaux2020adaptive}, a Type~II Hamiltonian Taylor Variational Integrator (HTVI) was derived following the strategy from Section~\ref{section: Variational Integrators} with $\rho = 0 $ for the adaptive approach $p \rightarrow \mathring{p}$-Bregman Hamiltonian,
\begin{equation} 
	\begin{aligned}
		\bar{H}_{p\rightarrow \mathring{p}}(\bar{q},\bar{r}) & =   \frac{p^2}{2\mathring{p}\mathfrak{q}^{p+\mathring{p}/p} }  \langle r , r \rangle  + \frac{Cp^2}{\mathring{p}} \mathfrak{q}^{2p-\mathring{p}/p} f(q) + \frac{p}{\mathring{p}} \mathfrak{r} \mathfrak{q}^{1-\mathring{p}/p}  .
	\end{aligned}
\end{equation}
This adaptive HTVI is the most natural Hamiltonian analogue of the LTVI described in Section~\ref{section: LTVI Algorithm}, and its updates are given by
	\begin{align}
		\mathfrak{q}_{k+1}  & = \mathfrak{q}_k + h \frac{p}{\mathring{p}} \mathfrak{q}_k^{1-\mathring{p}/p} , \\
	r_{k+1} & = r_k -  \frac{p^2}{\mathring{p}} Ch \mathfrak{q}_k^{2p-\mathring{p}/p} \nabla f(q_k), \\
	q_{k+1} &= q_k +  \frac{p^2}{\mathring{p} } h\mathfrak{q}_k^{-p-\mathring{p}/p} r_{k+1}. 
\end{align}
When $\mathring{p} = p$, it reduces to the direct HTVI,
	\begin{align}
	\mathfrak{q}_{k+1}  & = \mathfrak{q}_k + h , \\
	r_{k+1} & = r_k -   hCp \mathfrak{q}_k^{2p-1} \nabla f(q_k), \\
	q_{k+1} &= q_k +   hp\mathfrak{q}_k^{-p-1} r_{k+1}. 
\end{align} 
\hfill

\subsection{Numerical Experiments}

Numerical experiments using the numerical methods presented in the previous section have been conducted to minimize the quartic function,
\begin{equation}\label{Quartic}
	f(x) = \left[(x-1)^\top  \Sigma (x-1) \right]^2,
\end{equation} 
where $x\in \mathbb{R}^{d}$ and  $\Sigma_{ij} =0.9^{|i-j|}$. This convex function achieves its minimum value $0$ at $x^*=1$.   \\

\begin{figure}[!h]
	\centering
		\includegraphics[width=1\textwidth]{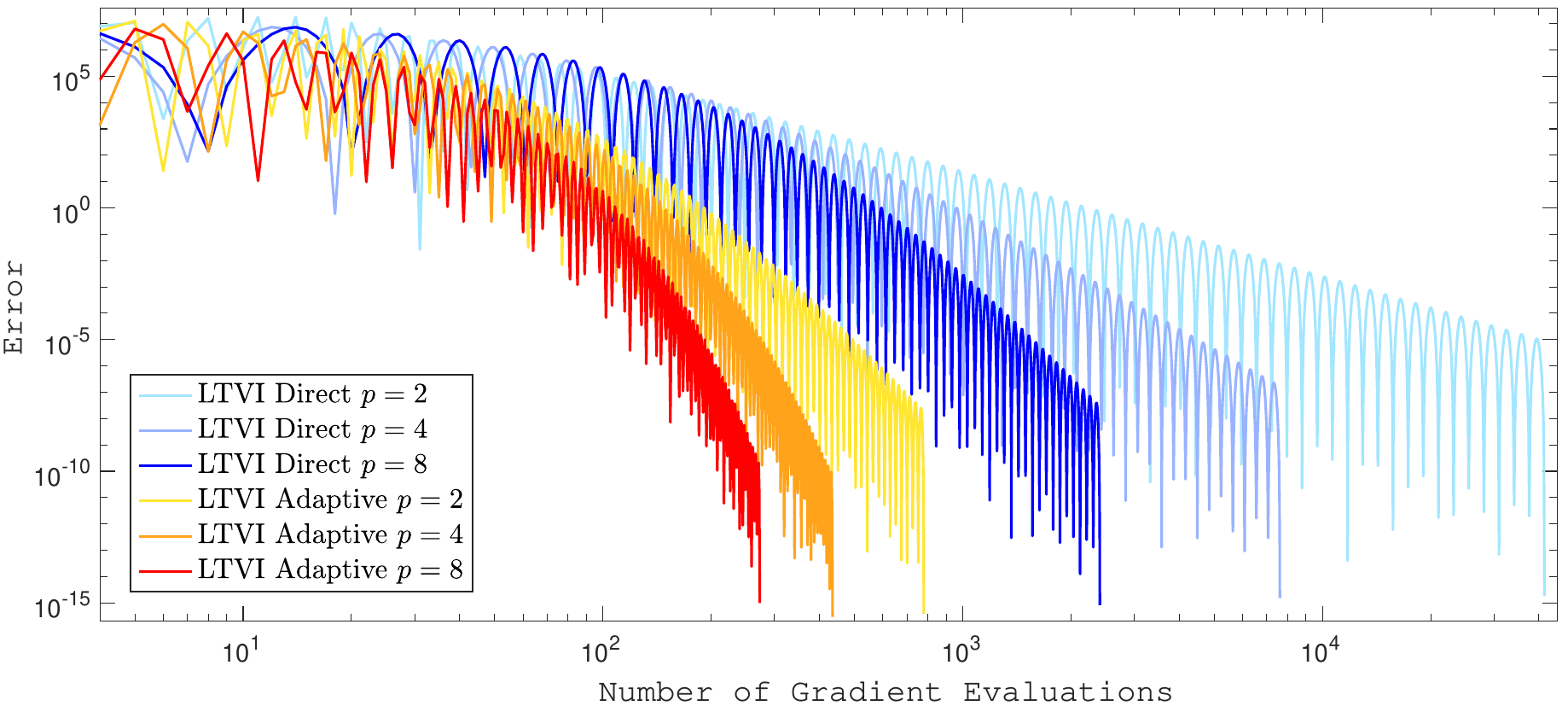}
	\caption{Comparison of the direct and adaptive approaches for the LTVI algorithm, when applied to the quartic function~\eqref{Quartic}. \label{fig:LTVI_ADvsDirect}}
\end{figure}

As was observed in~\cite{duruisseaux2020adaptive} for the HTVI algorithm, the numerical experiments showed that a carefully tuned adaptive approach algorithm enjoys a significantly better rate of convergence and requires a much smaller number of steps to achieve convergence than the direct approach, as can be seen in Figure~\ref{fig:LTVI_ADvsDirect} for the LTVI methods. Although the adaptive approach requires a smaller fictive time-step~$h$ than the direct approach, the physical time-steps resulting from $t = \tau^{p/\mathring{p}}$ in the adaptive approach grow rapidly to values larger than the constant physical time-step of the direct approach. The results of Figure~\ref{fig:LTVI_ADvsDirect} also show that the adaptive and direct LTVI methods become more and more efficient as $p$ is increased, which was also the case for the HTVI algorithm in~\cite{duruisseaux2020adaptive}. \\

\begin{figure}[!ht]
	\centering
	\includegraphics[width=1\textwidth]{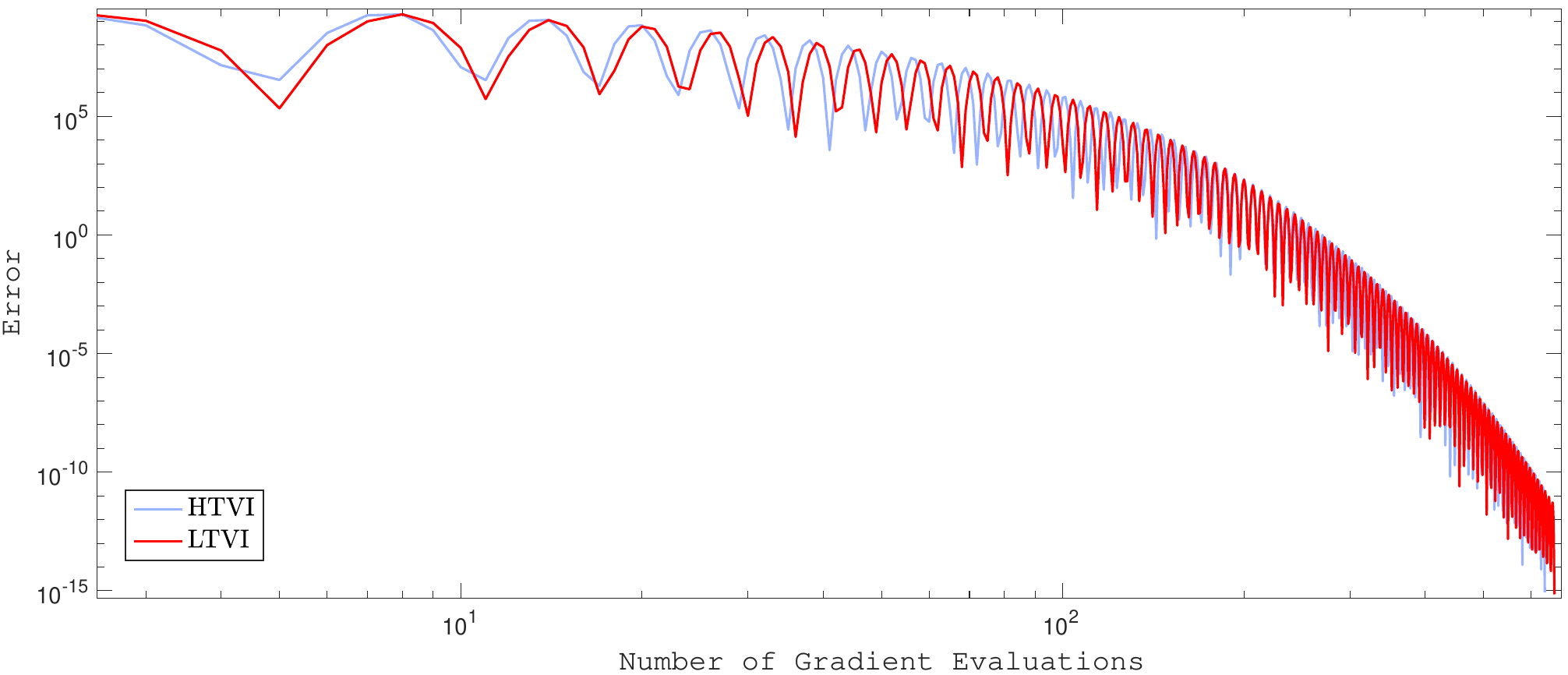}
	\caption{Comparison of the HTVI and LTVI algorithms with the same parameters. \label{fig: HTVI vs LTVI}} 
\end{figure}

The LTVI and HTVI algorithms presented in Section~\ref{section: Numerical Methods} perform empirically almost exactly in the same way for the same parameters and time-step, as can be seen for instance in Figure~\ref{fig: HTVI vs LTVI}. As a result, the computational analysis carried in~\cite{duruisseaux2020adaptive} for the HTVI algorithm extends to the LTVI algorithm. In particular, it was shown in~\cite{duruisseaux2020adaptive} that the HTVI algorithm is much more efficient than non-adaptive non-symplectic and adaptive non-symplectic integrators for the Bregman dynamics, and that it can be a competitive first-order explicit algorithm which can outperform certain popular optimization algorithms such as Nesterov's Accelerated Gradient~\cite{Nes83}, Trust Region Steepest Descent, ADAM~\cite{ADAM}, AdaGrad~\cite{AdaGrad}, and RMSprop~\cite{RMSprop}, for certain objective functions. Since the computational performance of the LTVI algorithm is almost exactly the same as that of the HTVI algorithm, this means that the LTVI algorithm is also much more efficient than non-symplectic integrators for the Bregman dynamics and can also be very competitive as a first-order explicit optimization algorithm. 

\section{Accelerated Optimization on More General Spaces} \label{section: Accelerated Optimization General Spaces}

\subsection{Motivation and Prior Work}

The variational framework for accelerated optimization on normed vector spaces from~\cite{WiWiJo16,duruisseaux2020adaptive} was extended to the Riemannian manifold setting in
\cite{Duruisseaux2022Riemannian} via a Riemannian $p$-Bregman Lagrangian  $\mathcal{L}_{p} :   T\mathcal{Q} \times \mathbb{R} \rightarrow  \mathbb{R} $ and a corresponding Riemannian $p$-Bregman Hamiltonian $\mathcal{H}_{p} :   T^*\mathcal{Q} \times \mathbb{R} \rightarrow  \mathbb{R} $, for $p>0$, of the form
\begin{equation} \label{BregmanLGeneral}
	\mathcal{L}_{p}(X,V,t) = \frac{t^{\frac{\zeta}{\lambda} p +1}}{2p} \langle V , V\rangle  - Cpt^{(\frac{\zeta}{\lambda} +1)p-1} f(X),  \end{equation} 
\begin{equation} \label{BregmanHGeneral}
	\mathcal{H}_{p}(X,R,t)= \frac{p}{2t^{\frac{\zeta}{\lambda}  p +1}} \llangle R , R\rrangle + Cpt^{(\frac{\zeta}{\lambda}  +1)p-1} f(X),
\end{equation}
where $\zeta$ and $\lambda$ are constants having to do with the curvature of the manifold and the convexity of the objective function $f$. These yield the associated $p$-Bregman Euler--Lagrange equations,
\begin{equation}  \nabla_{\dot{X}}\dot{X}  + \frac{\zeta p +\lambda}{\lambda t} \dot{X} + Cp^2t^{p-2} \text{grad}f(X) = 0. \end{equation} 
Here, $\text{grad}f$ denotes the Riemannian gradient of $f$, $\nabla_X Y $ is the covariant derivative of $Y$ along~$X$, and $  \llangle \cdot , \cdot \rrangle $ is the fiber metric on $T^* \mathcal{Q}$ induced by the Riemannian metric~$\langle \cdot , \cdot \rangle$ on $\mathcal{Q}$ whose local coordinates representation is the inverse of the local representation of $\langle \cdot , \cdot \rangle$. See~\cite{Absil2008,Boumal2020,MaRa1999,Lee2019,Jost2017,Duruisseaux2022Riemannian} for a more detailed description of these notions from Riemannian geometry and of this Riemannian variational framework for accelerated optimization. Note that some work was done on accelerated optimization via numerical integration of Bregman dynamics in the Lie group setting \cite{Lee2021,Tao2020} before the theory for more general Bregman families on Riemannian manifolds was established in~\cite{Duruisseaux2022Riemannian}.

 It was shown in~\cite{Duruisseaux2022Riemannian} that solutions to the  $p$-Bregman Euler--Lagrange equations converge to a minimizer of $f$ at a convergence rate of $\mathcal{O}(1/t^p)$, under suitable assumptions, and proven that time-rescaling a solution to the $p$-Bregman Euler--Lagrange equations via $\tau(t) =t^{\mathring{p}/p}$ yields a solution to the $\mathring{p}$-Bregman Euler--Lagrange equations. As a result, the adaptive approach involving the Poincar\'e transformation was extended to the Riemannian manifold setting via the adaptive approach Riemannian $p\rightarrow \mathring{p}$ Bregman Hamiltonian,
\begin{equation}  \label{H-R-Adaptive}
	\begin{aligned}
		\bar{\mathcal{H}}_{p \rightarrow \mathring{p}} \left(\bar{Q},\bar{R} \right)  & = 	 \frac{p^2}{2\mathring{p} \mathfrak{Q}^{\frac{\zeta}{\lambda} p +\frac{\mathring{p}}{p}}}  \llangle R , R\rrangle +  \frac{Cp^2}{\mathring{p}}\mathfrak{Q}^{\left(\frac{\zeta}{\lambda} +1\right)p-\frac{\mathring{p}}{p}} f(Q) + \frac{p}{\mathring{p}} \mathfrak{Q}^{1-\frac{\mathring{p}}{p}}   \mathfrak{R} .
	\end{aligned}
\end{equation}

This adaptive framework was exploited using discrete variational integrators incorporating holonomic constraints~\cite{Duruisseaux2022Constrained} and projection-based variational integrators~\cite{Duruisseaux2022Projection}. Both these strategies relied on embedding the Riemannian manifolds into an ambient Euclidean space. Although the Whitney and Nash Embedding Theorems~\cite{Whitney1944_2,Whitney1944_1,Nash1956} imply that there is no loss of generality when studying Riemannian manifolds only as submanifolds of Euclidean spaces, designing intrinsic methods that would exploit and preserve the symmetries and geometric properties of the Riemannian manifold and of the problem at hand could have advantages, both in terms of computation and in terms of improving our understanding of the acceleration phenomenon on Riemannian manifolds. 

Developing an intrinsic extension of Hamiltonian variational integrators to manifolds would require some additional work, since the current approach involves Type~II/III generating functions $H_d^+(q_k, p_{k+1})$, $H_d^-(p_k, q_{k+1})$, which depend on the position at one boundary point, and on the momentum at the other boundary point. However, this does not make intrinsic sense on a manifold, since one needs the base point in order to specify the corresponding cotangent space, and one should ideally consider a Hamiltonian variational integrator construction based on the discrete generalized energy in the context of discrete Dirac mechanics~\cite{LeOh2008,YoMa2006a,YoMa2006b}.  

On the other hand, Lagrangian variational integrators involve a Type~I generating function $L_d(q_k, q_{k+1})$ which only depends on the position at the boundary points and is therefore well-defined on manifolds, and many Lagrangian variational integrators have been derived on Riemannian manifolds, especially in the Lie group~\cite{LeLeMc2007a,LeLeMc2007b, BRMa2009,HaLe2012,Nordkvist2010,Le2004,Hussein2006,LeLeMc2005,LeeThesis} and homogeneous space~\cite{LeLeMc2009b} settings. The time-adaptive framework developed in this paper now makes it possible to design time-adaptive Lagrangian integrators for accelerated optimization on these more general spaces, where it is more natural and easier to work on the Lagrangian side than on the Hamiltonian side. \\

\subsection{Accelerated Optimization on Lie Groups}

Although it is possible to work on Riemannian manifolds, we will restrict ourselves to Lie groups for simplicity of exposition since there is more literature available on Lie group integrators than Riemannian integrators. Note as well that prior work is available on accelerated optimization via numerical integration of Bregman dynamics in the Lie group setting \cite{Lee2021,Tao2020}.

Here, we will work in the setting introduced in~\cite{Lee2021}. The setting of \cite{Tao2020} can be thought of as a special case of the more general Lie group framework for accelerated optimization presented here. Consider a $n$-dimensional Lie group $G$ with associated Lie algebra $\mathfrak{g} = T_eG$, and a left-trivialization of the tangent bundle of the group $TG \simeq  G \times  \mathfrak{g}$, via  $(q, \dot{q}) \mapsto  (q, \text{L}_{q^{-1}} \dot{q}) \equiv (q,\xi)$, where $\text{L} : G\times G \rightarrow G$ is the left action defined by $\text{L}_q h = qh$ for all $q,h\in G$. Suppose that $\mathfrak{g}$ is equipped with an inner product which induces an inner product on $T_qG$ via left-trivialization,
 $$  (v \bullet w )_{T_qG} =  ( \text{T}_q \text{L}_{q^{-1}} v \bullet  \text{T}_q \text{L}_{q^{-1}}w )_{\mathfrak{g}},  \qquad  \forall v, w \in T_qG. $$ 
 With this inner product, we identify $ \mathfrak{g} \simeq \mathfrak{g}^*$ and $T_q G \simeq T_q^* G \simeq G \times \mathfrak{g}^*$ via the Riesz representation. Let $\mathbf{J} : \mathfrak{g} \rightarrow \mathfrak{g}^*$ be chosen such that $( \mathbf{J}(\xi) \bullet \zeta ) $ is positive-definite and symmetric as a bilinear form of $\xi,\zeta \in \mathfrak{g}  $.  Then, the metric $\langle \cdot,\cdot  \rangle : \mathfrak{g} \times \mathfrak{g} \rightarrow \mathbb{R}$ with $\langle \xi , \zeta \rangle =   ( \mathbf{J}(\xi) \bullet  \zeta ) $ serves as a left-invariant Riemannian metric on $G$. The adjoint and $\text{ad}$ operators are denoted by $\text{Ad}_q : \mathfrak{g} \rightarrow \mathfrak{g}$ and $\text{ad}_{\xi} :\mathfrak{g} \rightarrow \mathfrak{g}$, respectively. We refer the reader to~\cite{MaRa1999,Lee2017,Gallier2020} for a more detailed description of Lie group theory and mechanics on Lie groups. 
 
 As mentioned earlier, there is a lot of literature available on Lie group integrators. We refer the reader to \cite{Iserles2000,Celledoni2014,Celledoni2022,Christiansen2011} for very thorough surveys of the literature on Lie group methods, which acknowledge all the foundational contributions leading to the current state of Lie group integrator theory. In particular, the Crouch and Grossman approach~\cite{Crouch1993}, the Lewis and Simo approach~\cite{Lewis1994}, Runge--Kutta--Munthe--Kaas methods \cite{MuntheKaas1998,MuntheKaas1999,MuntheKaas1995,Casas2003}, Magnus and Fer expansions~\cite{Zanna1999,Blanes2008,Iserles1999}, and commutator-free Lie group methods~\cite{Celledoni2003} are outlined in these surveys. Variational integrators have also been derived on the Lagrangian side in the Lie group setting~\cite{LeLeMc2007a,LeLeMc2007b, BRMa2009,HaLe2012,Nordkvist2010,Le2004,Hussein2006,LeLeMc2005,LeeThesis}.

We now introduce a discrete variational formulation of time-adaptive Lagrangian mechanics on Lie groups. Suppose we are given a partition $0 = \tau_0 < \tau_1 < \ldots < \tau_N = \mathcal{T}$ of the interval $[0,\mathcal{T}]$, and a discrete curve in $G \times \mathbb{R} \times \mathbb{R} $ denoted by $\{  (q_k, \mathfrak{q}_k, \lambda _k) \}_{k=0}^{N}$ such that $q_k \approx q(\tau_k)$,   $\mathfrak{q}_k \approx \mathfrak{q}(\tau_k)$, and $\lambda_k \approx \lambda(\tau_k)$. The discrete kinematics equation is chosen to be
		\begin{equation}
			q_{k+1} = q_k f_k,
		\end{equation}
where $f_k \in G$ represents the relative update over a single step.

 Consider the discrete action functional,
\begin{equation}
\bar{\mathfrak{S}}_d \left(\{  (q_k, \mathfrak{q}_k, \lambda _k) \}_{k=0}^{N} \right) =   \sum_{k=0}^{N-1}{  \left[  L_d(q_k ,f_k , \mathfrak{q}_k, \mathfrak{q}_{k +1} )   -  \lambda_{k} \frac{ \mathfrak{q}_{k+1} - \mathfrak{q}_k}{\tau_{k+1} - \tau_{k}}  + \lambda_k g(\mathfrak{q}_k)   \right]  \frac{ \mathfrak{q}_{k+1} - \mathfrak{q}_k}{\tau_{k+1} - \tau_{k}}    },
\end{equation}
	where,
\begin{equation}
	L_d(q_k ,f_k , \mathfrak{q}_k,  \mathfrak{q}_{k+1} )  \approx  \ext_{ \substack{ (q,\mathfrak{q} ) \in C^2 ([\tau_{k} , \tau_{k+1} ] ,G\times \mathbb{R} ) \\ (q, \mathfrak{q})(\tau_k) = (q_k, \mathfrak{q}_k)  , \text{  }  (q, \mathfrak{q})(\tau_{k+1}) = \left(q_{k} f_k, \mathfrak{q}_{k+1} \right) }  }   \text{  } \int_{\tau_k}^{\tau_{k+1}}{L\left(q, \frac{\xi }{g(\mathfrak{q})}  , \mathfrak{q} \right) d\tau }.
\end{equation}
\hfill \\

 		We can derive the following result which relates a discrete Type~I variational principle to a set of discrete Euler--Lagrange equations:

	\begin{theorem} \label{Theorem: discrete EL equations Lie Group} The Type~I discrete Hamilton's variational principle,
		\begin{equation}
			\delta \bar{\mathfrak{S}}_d \left(\{  (q_k, \mathfrak{q}_k, \lambda _k) \}_{k=0}^{N} \right)  = 0,
		\end{equation}
		where,
\begin{equation}
	\bar{\mathfrak{S}}_d \left(\{  (q_k, \mathfrak{q}_k, \lambda _k) \}_{k=0}^{N} \right) =   \sum_{k=0}^{N-1}{  \left[  L_d(q_k ,f_k , \mathfrak{q}_k, \mathfrak{q}_{k +1} )   -  \lambda_{k} \frac{ \mathfrak{q}_{k+1} - \mathfrak{q}_k}{\tau_{k+1} - \tau_{k}}  + \lambda_k g(\mathfrak{q}_k)   \right]  \frac{ \mathfrak{q}_{k+1} - \mathfrak{q}_k}{\tau_{k+1} - \tau_{k}}    },
\end{equation}
		is equivalent to the discrete extended Euler--Lagrange equations,
		\begin{equation}  \mathfrak{q}_{k+1} =  \mathfrak{q}_k + (\tau_{k+1} - \tau_{k} )g(\mathfrak{q}_k), \end{equation}  
\begin{equation}  \emph{Ad}^*_{f_k^{-1}} \left( \emph{T}_e^* \emph{L}_{f_k}   D_2 L_{d_k}   \right)  =   \emph{T}_e^* \emph{L}_{q_k}  D_1  L_{d_k}     + \frac{\tau_{k+1} - \tau_{k}}{ \mathfrak{q}_{k+1} - \mathfrak{q}_k}     \frac{ \mathfrak{q}_{k} - \mathfrak{q}_{k-1}}{\tau_{k} - \tau_{k-1}}    \emph{T}_e^* \emph{L}_{f_{k-1}}   D_2 L_{d_{k-1}}    ,    
\end{equation}  
\small 
\begin{equation} 
				\begin{aligned}
		&    \left[    D_3 L_{d_k} + \lambda_{k} \frac{ 1}{\tau_{k+1} - \tau_{k}}  + \lambda_k \nabla  g(\mathfrak{q}_k)  \right] \frac{ \mathfrak{q}_{k+1} - \mathfrak{q}_k}{\tau_{k+1} - \tau_{k}}      -  \frac{ 1}{\tau_{k+1} - \tau_{k}} \left[    L_{d_k} -  \lambda_{k} \frac{ \mathfrak{q}_{k+1} - \mathfrak{q}_k}{\tau_{k+1} - \tau_{k}}  + \lambda_k g(\mathfrak{q}_k)      \right]  \\ & \qquad    + \left[    D_4L_{d_k}  - \lambda_{k-1} \frac{ 1}{\tau_{k} - \tau_{k-1}}   \right] \frac{ \mathfrak{q}_{k} - \mathfrak{q}_{k-1}}{\tau_{k} - \tau_{k-1}}   +  \frac{ 1}{\tau_{k} - \tau_{k-1}} \left[    L_{d_{k-1}}  -  \lambda_{k-1} \frac{ \mathfrak{q}_{k} - \mathfrak{q}_{k-1}}{\tau_{k} - \tau_{k-1}}  + \lambda_{k-1} g(\mathfrak{q}_{k-1})      \right]     =0,
	\end{aligned} \end{equation}  
where $L_{d_k} $ denotes $ L_d(q_k ,f_k , \mathfrak{q}_k, \mathfrak{q}_{k +1} )   $.
\normalsize
		\proof{See Appendix \ref{Appendix: Discrete Proof Lie Group}.		\qed \\}
	\end{theorem}  

Now, define 
\begin{align} \mathfrak{p}_k & = - D_3 L_d(q_k ,f_k , \mathfrak{q}_k, \mathfrak{q}_{k +1} ) \end{align} and \begin{align} \mu_{k} & =  \text{Ad}^*_{f_k^{-1}} \left( \text{T}_e^* \text{L}_{f_k}   D_2 L_d(q_k ,f_k , \mathfrak{q}_k, \mathfrak{q}_{k +1} )   \right) -  \text{T}_e^* \text{L}_{q_k}  D_1  L_d(q_k ,f_k , \mathfrak{q}_k, \mathfrak{q}_{k +1} ) . \end{align} 
Then,
\begin{align}
	\mu_{k+1} & = \frac{\tau_{k+2} - \tau_{k+1}}{ \mathfrak{q}_{k+2} - \mathfrak{q}_{k+1}}     \frac{ \mathfrak{q}_{k+1} - \mathfrak{q}_{k}}{\tau_{k+1} - \tau_{k}}    \text{T}_e^* \text{L}_{f_{k}}   D_2 L_d(q_k ,f_k , \mathfrak{q}_k, \mathfrak{q}_{k +1} )   .\end{align}  

\noindent With these definitions, if we use a constant time-step $h$ in fictive time $\tau$ and substitute $g(\mathfrak{q}) = \frac{p}{\mathring{p}} \mathfrak{q}^{1-\mathring{p}/p}  $, the discrete Euler--Lagrange equations can be rewritten as
\begin{align}  \mu_{k} & =  \text{Ad}^*_{f_k^{-1}} \left( \text{T}_e^* \text{L}_{f_k}   D_2 L_d(q_k ,f_k , \mathfrak{q}_k, \mathfrak{q}_{k +1} )   \right) -  \text{T}_e^* \text{L}_{q_k}  D_1  L_d(q_k ,f_k , \mathfrak{q}_k, \mathfrak{q}_{k +1} ) ,  \\  
	\mu_{k+1} & =  \frac{  \mathfrak{q}_k^{1-\mathring{p}/p}   }{  \mathfrak{q}_{k+1}^{1-\mathring{p}/p}   }  \text{T}_e^* \text{L}_{f_{k}}   D_2 L_d(q_k ,f_k , \mathfrak{q}_k, \mathfrak{q}_{k +1} ),   \\ \mathfrak{q}_{k+1} & =  \mathfrak{q}_k +h \frac{p}{\mathring{p}} \mathfrak{q}_k^{1-\mathring{p}/p}  , \\
	\mathfrak{p}_{k+1} & = \frac{\mathring{p}\left[ \lambda_{k+1} - \lambda_k  + L_{d_k} -  L_{d_{k+1}}    \right] }{h p \mathfrak{q}_{k+1}^{1-\mathring{p}/p} }   +   	  D_4 L_{d_k} +  \frac{\lambda_{k+1}}{  \mathfrak{q}_{k+1}} \left(1-\frac{\mathring{p}}{p} \right)  	   .  \end{align} 
\hfill 
   
   \normalsize
   	In the Lie group setting, the Riemannian $p$-Bregman Lagrangian becomes
   \begin{equation} \label{eq: Lie Group Bregman Lagrangian}
   	\mathcal{L}_{p}(q,\xi,t) = \frac{t^{ \kappa  p +1}}{2p} \langle \xi , \xi \rangle  - Cpt^{(\kappa  +1)p-1} f(q),  \end{equation} 
   with corresponding Euler--Lagrange equation, 
   \begin{equation}
   	\frac{d \mathbf{J}(\xi)}{dt} + \frac{\kappa p +1}{t} \mathbf{J}(\xi)  - \text{ad}_{\xi}^{*} \mathbf{J}(\xi) + Cp^2 t^{p-2} \nabla_{\text{L}} f(q) = 0,
   \end{equation} 
   where $\nabla_{\text{L}} f$ is the left-trivialized derivative of $f$, given by $
   \nabla_{\text{L}} f(q) = \text{T}_e^* \text{L}_q(\text{D}_q f(q))$.  We then consider the discrete Lagrangian,
   \begin{align}
   	L_d(q_k ,f_k , \mathfrak{q}_k,  \mathfrak{q}_{k+1} )  &  =  \frac{\mathfrak{q}_k^{\kappa p +1}}{hp (g(\mathfrak{q}_k))^2}  T_d (f_k)  - Chp\mathfrak{q}_k^{(\kappa +1)p-1} f(q_k),
   \end{align} 
   where $T_d(f_k) \approx \frac{1}{2} \langle h\xi_k , h\xi_k \rangle$, which approximates
   \begin{equation}
   	L_d(q_k ,f_k , \mathfrak{q}_k,  \mathfrak{q}_{k+1} )  \approx  \ext_{ \substack{ (q,\mathfrak{q} ) \in C^2 ([\tau_{k} , \tau_{k+1} ] ,G\times \mathbb{R} ) \\ (q, \mathfrak{q})(\tau_k) = (q_k, \mathfrak{q}_k)   , \text{  }  (q, \mathfrak{q})(\tau_{k+1}) = (q_{k} f_k, \mathfrak{q}_{k+1}) }  }   \text{  } \int_{\tau_k}^{\tau_{k+1}}{L\left(q, \frac{\xi }{g(\mathfrak{q})}  , \mathfrak{q} \right) d\tau }.
   \end{equation}

\subsection{Numerical Experiment on $\text{SO}(3)$}

We work on the 3-dimensional Special Orthogonal group,
\begin{equation} \text{SO}(3) = \{ R \in \mathbb{R}^{3\times 3} |  R^\top R = I_{3\times 3} , \det{(R)} = 1 \}.
\end{equation}
Its Lie algebra is
\begin{equation} \mathfrak{so}(3) = \{ S \in \mathbb{R}^{3\times 3} |  S^\top  = -S \},
\end{equation}
with the matrix commutator as the Lie bracket. We have an identification between $\mathbb{R}^{3} $ and $ \mathfrak{so}(3)$ given by the hat map $\hat{\cdot} : \mathbb{R}^{3} \rightarrow \mathfrak{so}(3) $, defined such that $\hat{x} y = x\times y$ for any $x,y\in \mathbb{R}^3$. The inverse of the hat map is the vee map $(\cdot)^{\vee} :  \mathfrak{so}(3) \rightarrow \mathbb{R}^{3}$. The inner product on $\mathfrak{so}(3) $ is given by
\begin{equation}
	\left(  \hat{\eta} \bullet  \hat{\xi}  \right)_{\mathfrak{so}(3)} = \frac{1}{2} \text{Trace} \left(\hat{\eta}^\top \hat{\xi} \right) = \eta^\top \xi,
\end{equation}
and the metric is chosen so that 
\begin{equation}
	\langle \hat{\eta} , \hat{\xi}   \rangle = \left( \mathbf{J}( \hat{\eta})  \bullet  \hat{\xi}  \right)_{\mathfrak{so}(3)} = \text{Trace} \left(\hat{\eta}^\top J_d \hat{\xi} \right) = \eta^\top J \xi ,
\end{equation}
where $J \in \mathbb{R}^{3\times 3}$ is a symmetric positive-definite matrix and $J_d = \frac{1}{2} \text{Trace}(J) I_{3\times 3} - J$. \\
 On $\text{SO}(3)$, for any $u,v \in \mathbb{R}^3$ and $F \in \text{SO}(3)$,
\begin{equation}
	\text{ad}_{\hat{u}} \hat{v} = [\hat{u},\hat{v}] =  \hat{u} \hat{v} - \hat{v} \hat{u} = \widehat{ u \times v}, \qquad 
	\text{Ad}_F \hat{u} = F \hat{u} F^\top  = \widehat{Fu}.
\end{equation} 
Identifying $\mathfrak{so}(3)^* \simeq \mathfrak{so}(3) \simeq \mathbb{R}^3$, we have for any $u,v \in \mathbb{R}^3$ and $F \in \text{SO}(3)$ that
\begin{equation}
	\text{ad}_u v = \hat{u} v = u \times v,  \qquad    \text{ad}^*_u v  = -\hat{u} v = v \times u,  \qquad 
	\text{Ad}_F u = Fu ,  \qquad 
	\text{Ad}_F^* u  = F^\top u.
\end{equation}   \\
\indent On $\text{SO}(3)$, the Riemannian $p$-Bregman Lagrangian becomes
\begin{equation}
	\mathcal{L}_{p}(R,\Omega,t) = \frac{t^{   p +1}}{2p} \Omega^\top  J\Omega  - Cpt^{2p-1} f(R), 
\end{equation}
and the corresponding Euler--Lagrange equations are given by
\begin{equation}
	J\dot{\Omega} + \frac{p +1}{t} J \Omega + \hat{\Omega} J \Omega + Cp^2 t^{p-2} \nabla_{\text{L}} f(R) = 0, \qquad  
	\dot{R} = R \hat{\Omega}.
\end{equation}  
The discrete kinematics equation is written as
\begin{equation}
	R_{k+1} = R_k F_k ,
\end{equation}
where $F_k \in \text{SO}(3)$, and $\kappa = 1$ so we get the discrete Lagrangian,
\begin{align}
	L_d(R_k ,F_k , \mathfrak{R}_k,  \mathfrak{R}_{k+1} )  &  =  \frac{ \mathring{p}^2}{h p^3 }  \mathfrak{R}_k^{ p -1 + 2\mathring{p}/p}  T_d(F_k)   - Chp\mathfrak{R}_k^{2p-1} f(R_k).
\end{align} 
As in~\cite{LeLeMc2007a,Lee2021}, the angular velocity is approximated by $\hat{\Omega}_k \approx \frac{1}{h} R_k^\top (R_{k+1}-R_k) = \frac{1}{h}  (F_{k}-I_{3\times 3}) $ so we can take
\begin{equation}
	T_d(F_k) = \text{Trace}\left( [I_{3\times 3} - F_k] J_d \right).
\end{equation} 
Differentiating this equation and using the identity $\text{Trace}(- \hat{x} A) = (A - A^\top )^{\vee} \cdot x$ yields
\begin{equation}
	\text{T}_I^* \text{L}_{F_k}\left(D_{F_k} T_d(F_k) \right) = \left( J_d F_k - F_k^\top J_d  \right)^{\vee}.
\end{equation}
Then, the discrete Euler--Lagrange equations for $\mu_k$ and $\mu_{k+1}$ become
\begin{align} \label{eq: Lie Group mu_k temp equation} \mu_{k} & = \frac{ \mathring{p}^2}{h p^3 }  \mathfrak{R}_k^{ p -1 + 2\mathring{p}/p}  \left( F_k J_d  - J_d F_k^\top   \right)^{\vee}  + C h p \mathfrak{R}_k^{2p-1} \nabla_{\text{L}} f(R_k),   \\
	\mu_{k+1} & = \frac{  \mathfrak{q}_k^{1-\mathring{p}/p}   }{  \mathfrak{q}_{k+1}^{1-\mathring{p}/p}   }  F_k^\top  \left[   \mu_k    -  C h p \mathfrak{R}_k^{2p-1}   \nabla_{\text{L}} f(R_k)     \right]  .\end{align} 
Now, equation \eqref{eq: Lie Group mu_k temp equation} can be solved explicitly when $J=I_{3\times 3}$ as described in~\cite{Lee2021}:
\begin{align}
	F_k = \exp{ \left(  \frac{\sin^{-1}{\| a_k \|}}{\|a_k \|} \hat{a}_k  \right)}, \qquad  \text{where} \quad  a_k = \frac{h p^3 } { \mathring{p}^2} \mathfrak{R}_k^{ 1-p - 2\mathring{p}/p} \left[   \mu_k -  C h p \mathfrak{R}_k^{2p-1} \nabla_{\text{L}} f(R_k) \right].
\end{align}
Therefore, we get the \textbf{Adaptive LLGVI} (Adaptive Lagrangian Lie Group Variational Integrator) 
\begin{align} 
	&F_k  = \exp{ \left(  \frac{\sin^{-1}{\| a_k \|}}{\|a_k \|} \hat{a}_k  \right)}, \text{  }  \text{ where }\text{ }  a_k = \frac{h p^3 } { \mathring{p}^2} \mathfrak{R}_k^{ 1-p - 2\mathring{p}/p} \left[   \mu_k -  C h p \mathfrak{R}_k^{2p-1} \nabla_{\text{L}} f(R_k) \right], \\
	& \mathfrak{R}_{k+1} =  \mathfrak{R}_k + h \frac{p}{\mathring{p}} \mathfrak{R}_k^{1-\mathring{p}/p} , \\
	&\mu_{k+1}  = \frac{  \mathfrak{R}_k^{1-\mathring{p}/p}   }{  \mathfrak{R}_{k+1}^{1-\mathring{p}/p}   }  F_k^\top  \left[   \mu_k    -  C h p \mathfrak{R}_k^{2p-1}   \nabla_{\text{L}} f(R_k)     \right] , \\ &R_{k+1}  = R_k F_k.
\end{align} 

We will use this integrator to solve the problem of minimizing the objective function,
\begin{equation} \label{eq: Wahba Objective Function}
	f(R) = \frac{1}{2} \| A-R \|_F^2 = \frac{1}{2} \left( \| A \|_F^2 + 3 \right) - \text{Trace}(A^\top R),
\end{equation}
over $R\in \text{SO}(3)$, where $\| \cdot \|_F$ denotes the Frobenius norm. Its left-trivialized gradient is given by 
\begin{equation}
	\nabla_{\text{L}} f(R)  = \left( A^\top R - R^\top A \right)^{\vee}.
\end{equation}
Minimizing this objective function appears in the least-squares estimation of attitude, commonly known as Wahba’s problem~\cite{Wahba1965}. The optimal attitude is explicitly given by 
\begin{equation} R^* = U \text{diag}\left[1, 1, \det(U V ) \right]V^\top , \end{equation} where $A = USV^\top $ is the singular value decomposition of $A$ with $U,V\in \text{O}(3)$ and $S$ diagonal. \\

\begin{figure}[!ht] 
	\centering
	\begin{minipage}[b]{0.48\textwidth}
		\includegraphics[width=\textwidth]{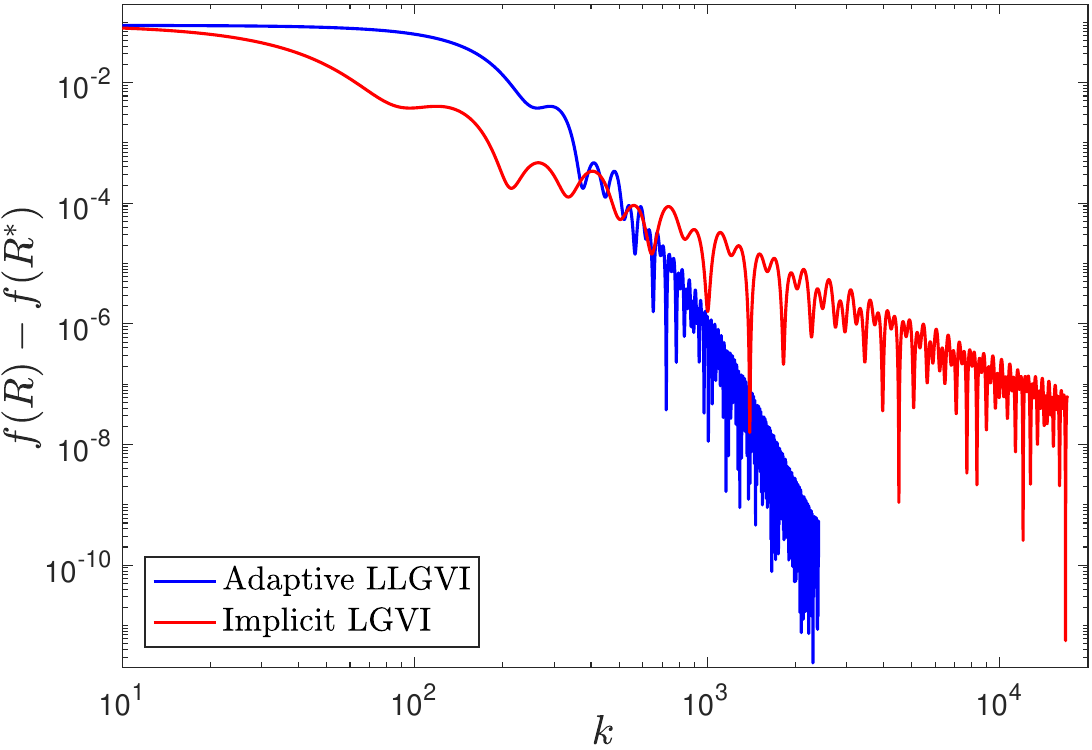}
	\end{minipage}
	\begin{minipage}[b]{0.48\textwidth}
		\includegraphics[width=\textwidth]{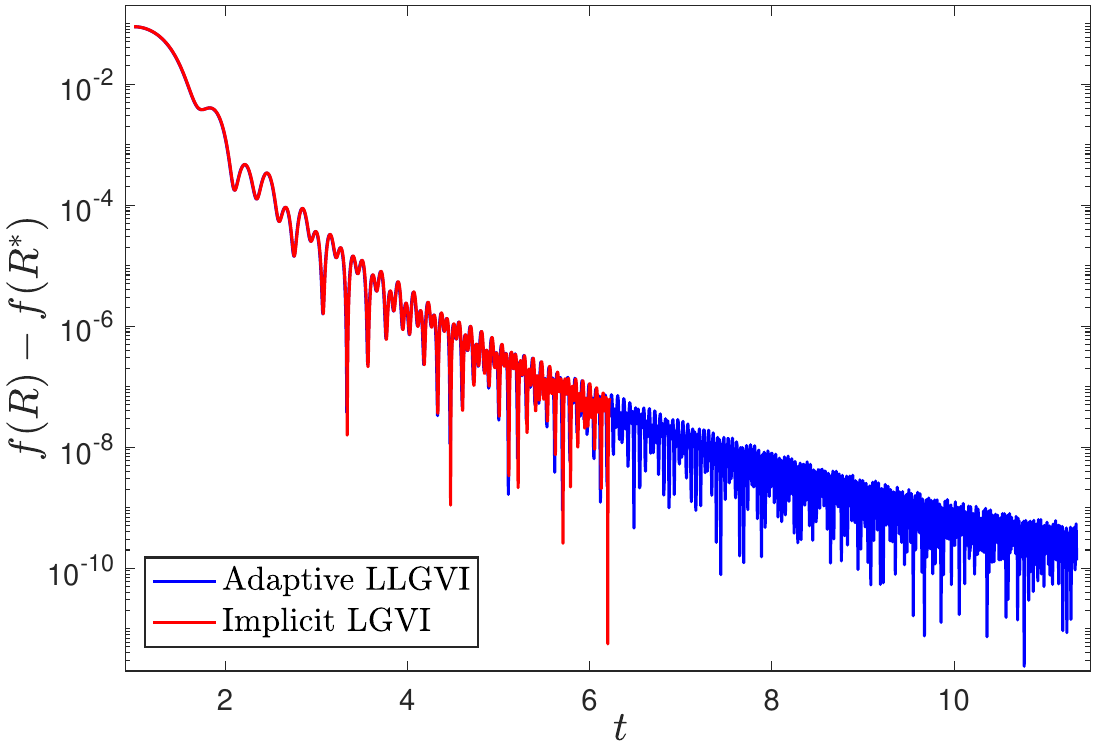}
	\end{minipage} \\
	\begin{minipage}[b]{0.48\textwidth}
		\includegraphics[width=\textwidth]{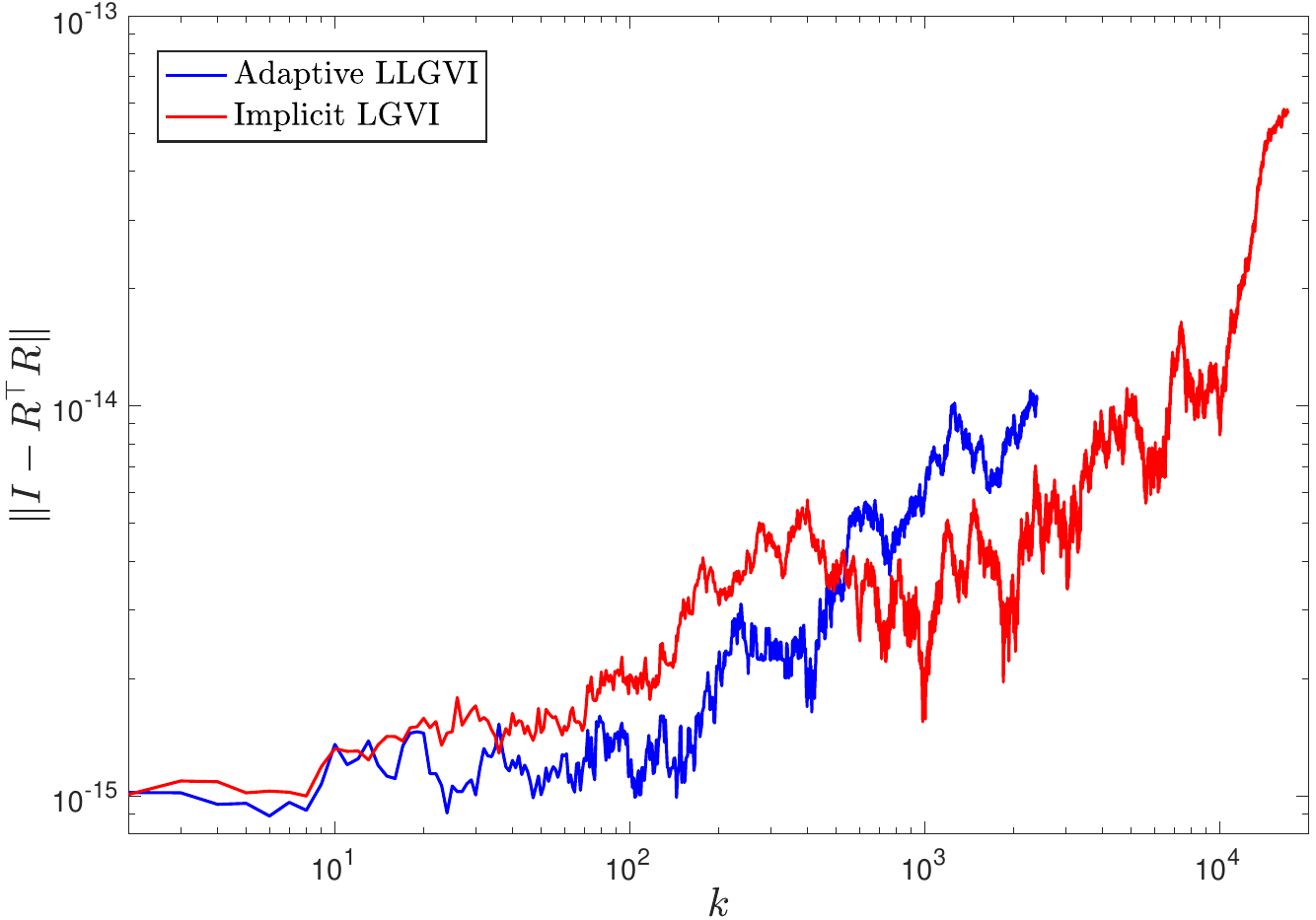}
	\end{minipage}
	\begin{minipage}[b]{0.48\textwidth}
		\includegraphics[width=\textwidth]{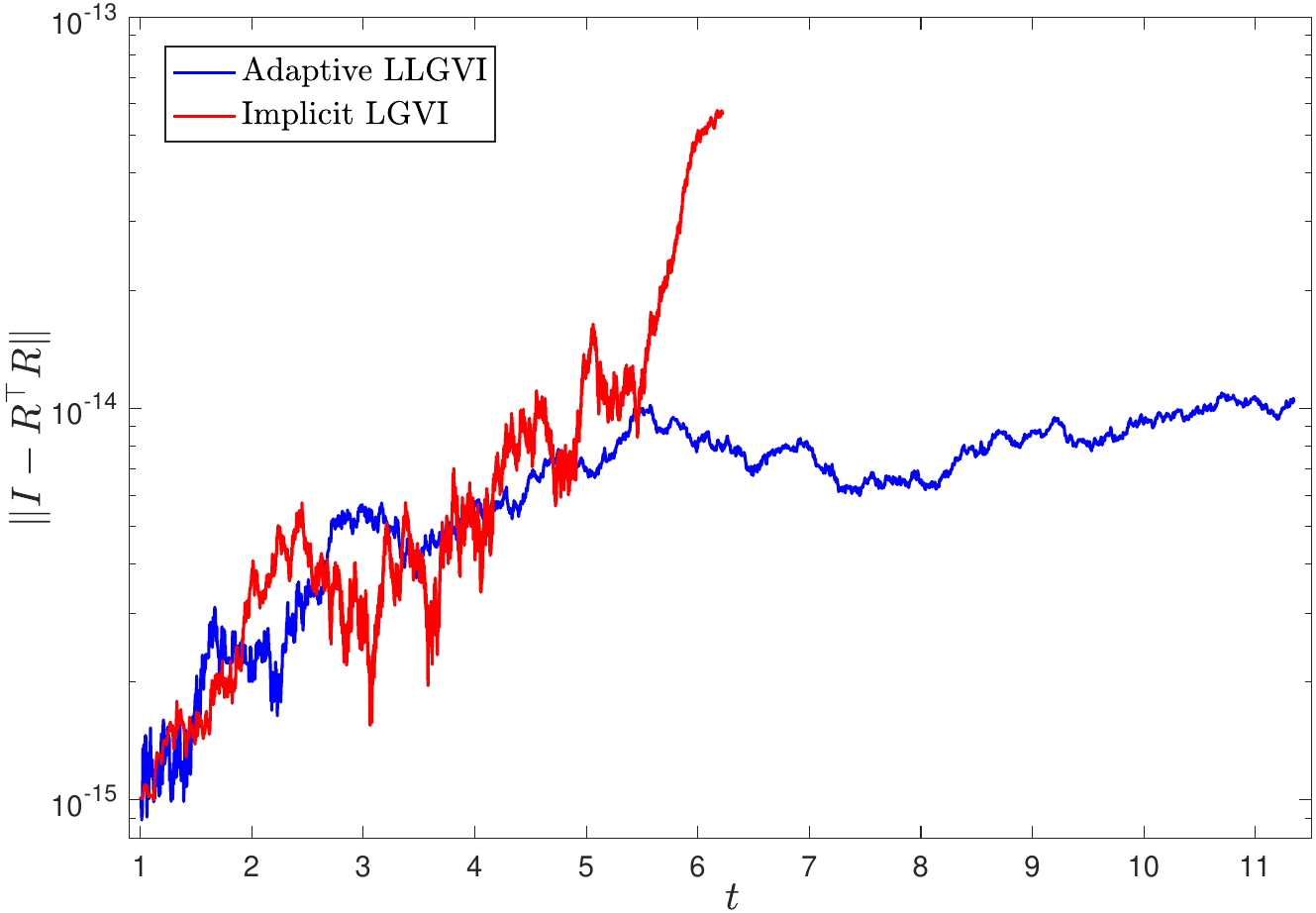}
	\end{minipage} \\
	\begin{minipage}[b]{0.49\textwidth}
		\includegraphics[width=\textwidth]{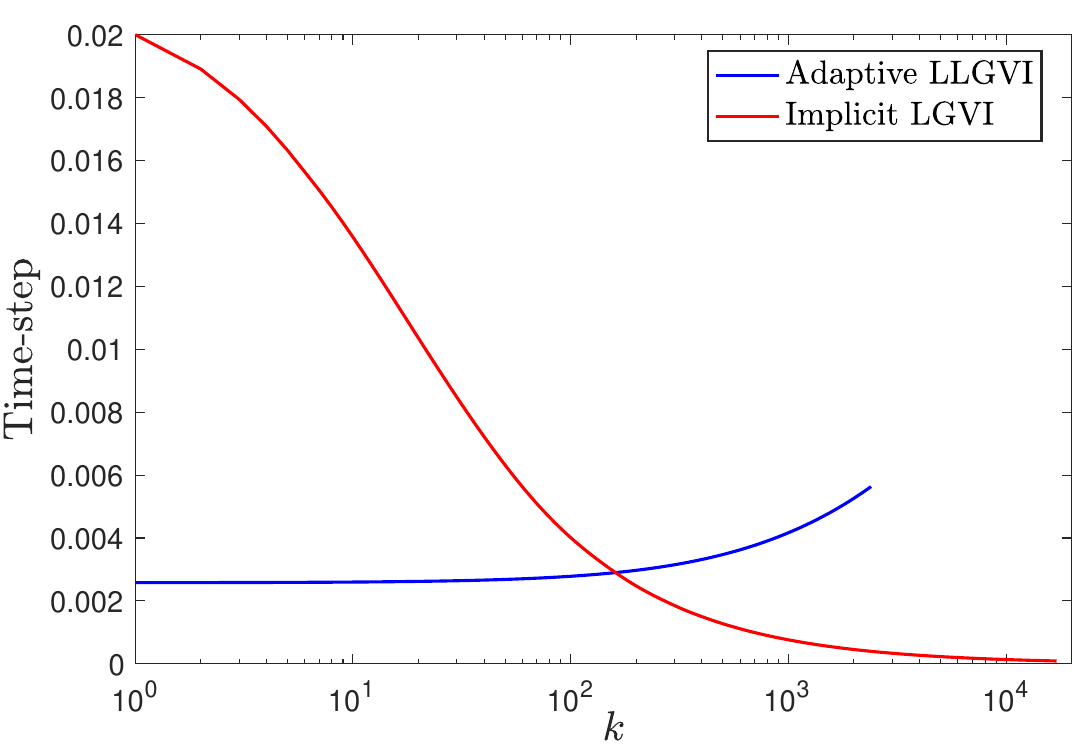}
	\end{minipage}
	\begin{minipage}[b]{0.48\textwidth}
		\includegraphics[width=\textwidth]{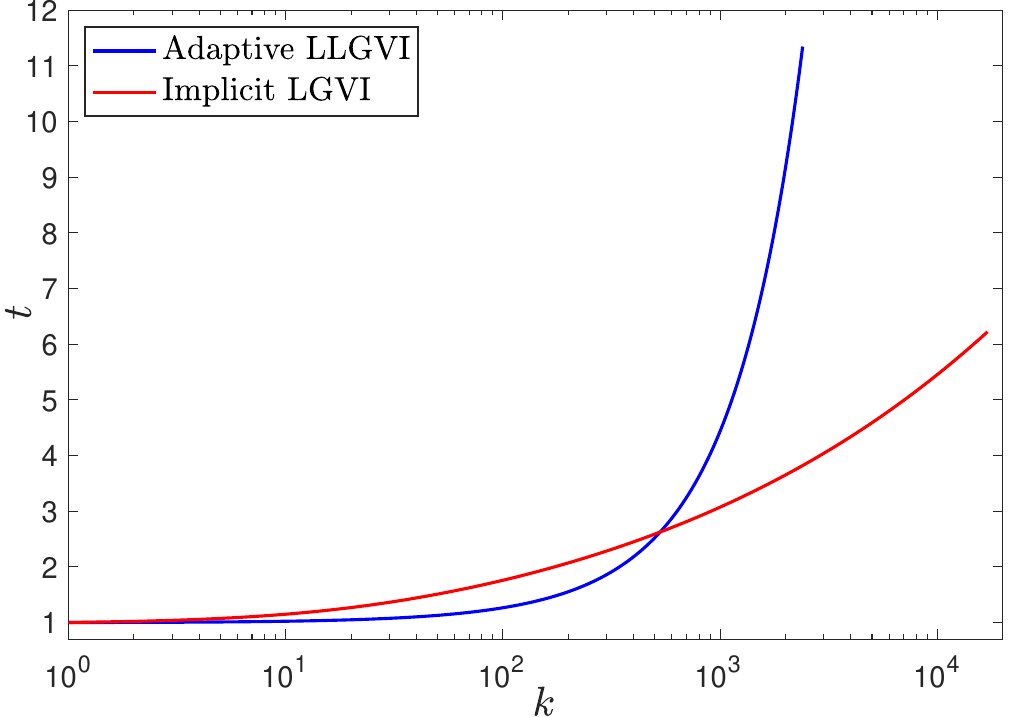}
	\end{minipage} 
\vspace*{-2mm}
	\caption{Comparison of the Adaptive LLGVI algorithm and of the Implicit LGVI algorithm from~\cite{Lee2021} with $p=6$, to solve Wahba's problem \eqref{eq: Wahba Objective Function}.}\label{fig: Lie Group Results}
\end{figure}

\begin{figure}[!ht] 
		\centering
		\includegraphics[width=0.99\textwidth]{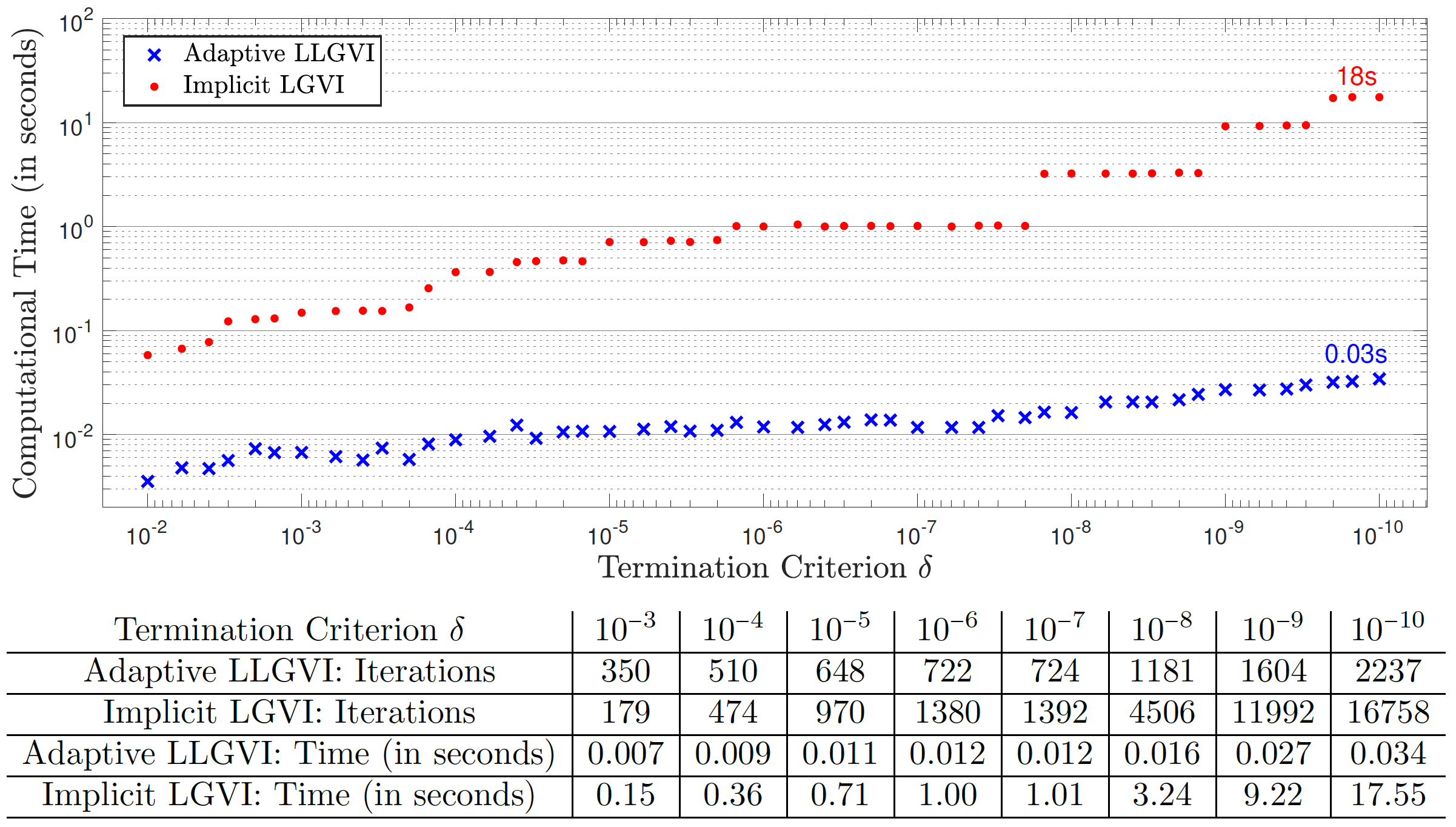} \caption{Time and number of iterations needed by the Adaptive LLGVI and Implicit LGVI algorithms with $p=6$, to satisfy the termination criterion \eqref{eq: Termination Criterion Lie Group} on Wahba's problem \eqref{eq: Wahba Objective Function}. }\label{fig: Computational Time}
	\end{figure}

We have tested the Adaptive LLGVI integrator on Wahba's problem against the Implicit Lie Group Variational Integrator (\textbf{Implicit LGVI}) from~\cite{Lee2021}. The Implicit LGVI is a Lagrangian Lie group variational integrator which adaptively adjusts the step size at every step. It should be noted that these two adaptive approaches use adaptivity in two fundamentally different ways: our Adaptive LLGVI method uses \textit{a priori} adaptivity based on known global properties of the family of differential equations considered (i.e. the time-rescaling symmetry of the family of Bregman dynamics), while the implicit method from~\cite{Lee2021} adapts the time-steps in an \textit{a posteriori} way, by solving a system of nonlinear equations coming from an extended variational principle. The results of our numerical experiments are presented in Figures~\ref{fig: Lie Group Results} and \ref{fig: Computational Time}. In these numerical experiments, we have used the termination criteria
\begin{equation} \label{eq: Termination Criterion Lie Group}
	| f(R_k)  - f(R^*) | < \delta  \qquad \text{and} \qquad |   f(R_k)  - f(R_{k-1}) | < \delta.
\end{equation}

We can see from Figure~\ref{fig: Lie Group Results} that both algorithms preserve the orthogonality condition $R_k^\top R_k = I_{3\times 3}$ very well. Now, we can observe from Figure~\ref{fig: Lie Group Results} that although both algorithms follow the same curve in time $t$, they do not travel along this curve at the same speed. Despite the fact that the Adaptive LLGVI algorithm initially takes smaller time-steps, those time-steps eventually become much larger than for the Implicit LGVI algorithm, and as a result, the Adaptive LLGVI algorithm achieves the termination criteria in a smaller number of iterations, which can also be seen more explicitly in the table from Figure~\ref{fig: Computational Time}. Unlike the Implicit LGVI algorithm, the Adaptive LLGVI algorithm is explicit, so each iteration is much cheaper and is therefore significantly faster, as can be seen from the running times displayed in Figure~\ref{fig: Computational Time}. Furthermore, the Adaptive LLGVI algorithm is significantly easier to implement.  \\

%Note however that the time-step in the Adaptive LLGVI algorithm increases indefinitely according to $ \mathfrak{R}_{k+1} =  \mathfrak{R}_k + h \frac{p}{\mathring{p}} \mathfrak{R}_k^{1-\mathring{p}/p} $. This can eventually lead to instability due to finite numerical precision, but there are relatively easy and cheap ways to fix this issue.

%
%
%\newpage 
%\begin{figure}[!ht] 
%	\centering
%	\includegraphics[width=0.99\textwidth]{ComputationalTime}  \vspace{-5mm}
%\end{figure}
%\begin{table}[!ht]
%	\begin{tabular}{c|c|c|c|c|c|c|c|c}
%		Termination Criterion
%		$\delta$                                      & $10^{-3} $ & $10^{-4}$ & $10^{-5}$ & $10^{-6}$ & $10^{-7}$ & $10^{-8}$ & $10^{-9}$ & $10^{-10}$ \\ \hline
%		Adaptive LLGVI: Iterations  &   350        &      510    &    648      &   722       &    724     &     1181      &    1604      &       2237     \\ \hline
%		Implicit LGVI: Iterations  &   179         &    474      &     970     &     1380      &    1392       &    4506       &   11992       &   16758         \\ \hline
%		Adaptive LLGVI: Time (in seconds)   & 0.007          &   0.009        &       0.011    &      0.012     &       0.012    &    0.016       &  0.027          &     0.034       \\ \hline
%		Implicit LGVI: Time (in seconds)      &        0.15    &    0.36      &   0.71        &    1.00       &    1.01      &  3.24         &   9.22        &  17.55       \\ \hline     
%	\end{tabular}
%\end{table}
%\newpage  

\section{Conclusion}

	A variational framework for accelerated optimization on vector spaces was introduced~\cite{WiWiJo16} by considering a family of time-dependent Bregman Lagrangian and Hamiltonian systems which is closed under time-rescaling. This variational framework was exploited in~\cite{duruisseaux2020adaptive} by using time-adaptive geometric Hamiltonian integrators to design efficient, explicit algorithms for symplectic accelerated optimization. It was observed that a careful use of adaptivity and symplecticity, which was possible on the Hamiltonian side thanks to the Poincar\'e transformation, could result in a significant gain in computational efficiency, by simulating higher-order Bregman dynamics using the computationally efficient lower-order Bregman integrators applied to the time-rescaled dynamics.

	The variational framework and time-adaptive approach on the Hamiltonian side were later extended to the Riemannian manifolds setting in~\cite{Duruisseaux2022Riemannian}. However, the current formulations of Hamiltonian variational integrators do not make sense intrinsically on manifolds, so this framework was only exploited using methods which take advantage of the structure of the Euclidean spaces in which the Riemannian manifolds are embedded~\cite{Duruisseaux2022Constrained,Duruisseaux2022Projection} instead of the structure of the Riemannian manifolds themselves. On the other hand, existing formulations of Lagrangian variational integrators are well-defined on manifolds, and many Lagrangian variational integrators have been derived on Riemannian manifolds, especially in the Lie group setting. This motivated exploring whether it is possible to construct a mechanism on the Lagrangian side which mimics the Poincar\'e transformation, since it is more natural and easier to work on the Lagrangian side on curved manifolds. 
	
	The usual correspondence between Hamiltonian and Lagrangian dynamics could not be exploited here since the Poincar\'e Hamiltonian is degenerate and therefore does not have a corresponding Lagrangian formulation. Instead, we introduced a novel derivation of the Poincar\'e transformation from a variational principle which gave us additional insight into the transformation mechanism and provided natural candidates for a time-adaptive framework on the Lagrangian side. Based on these observations, we constructed a theory of time-adaptive Lagrangian mechanics both in continuous and discrete time, and applied the resulting time-adaptive Lagrangian variational integrators to solve optimization problems by simulating Bregman dynamics, within the variational framework introduced in~\cite{WiWiJo16}. We observed empirically that our time-adaptive Lagrangian variational integrators performed almost exactly in the same way as the time-adaptive Hamiltonian variational integrators coming from the Poincar\'e framework of~\cite{duruisseaux2020adaptive}, whenever they are used with the same parameters and time-step. As a result, the computational analysis carried in~\cite{duruisseaux2020adaptive} for the HTVI algorithm extends to the LTVI algorithm, and thus the LTVI algorithm is much more efficient than non-symplectic integrators for the Bregman dynamics and can be a competitive first-order explicit algorithm since it can outperform commonly used optimization algorithms for certain objective functions.

 Finally, we showed that our time-adaptive Lagrangian approach extends naturally to more general spaces such as Riemannian manifolds and Lie groups without having to face the difficulties experienced on the Hamiltonian side, and we then applied time-adaptive Lie group Lagrangian variational integrators to solve an optimization problem on the three-dimensional Special Orthogonal group $\text{SO}(3)$. In particular, the resulting Lie group accelerated optimization algorithms were significantly faster and easier to implement than other recently proposed time-adaptive Lie group variational integrators for accelerated optimization.
 
 In future work, we will explore the issue of time-adaptive Lagrangian mechanics for more general monitor functions, using the primal-dual framework of Dirac mechanics. We will also study the convergence properties of the discrete-time algorithms, and try to better understand how to reconcile the Nesterov barrier theorem with the convergence properties of the continuous Bregman flows. It would also be useful to study the extent to which the practical considerations recently presented in~\cite{Duruisseaux2022Practical}, which significantly improved the computational performance of the symplectic optimization algorithms in the normed vector space setting, extend to the Riemannian manifold and Lie group settings with the Lagrangian Riemannian and Lie group variational integrators. \\

\section*{Acknowledgments} 
The authors would like to thank Taeyoung Lee and Molei Tao for sharing the code for the implicit Lie group variational integrator from \cite{Lee2021}, which was used in the numerical comparison with the method introduced in this paper. The authors would also like to thank the referees for their helpful comments and suggestions.

\section*{Funding}
The authors were supported in part by NSF under grants DMS-1345013, DMS-1813635, CCF-2112665, by AFOSR under grant FA9550-18-1-0288, and by the DoD under grant HQ00342010023 (Newton Award for Transformative Ideas during the COVID-19 Pandemic). \\

\appendix

\section{Proofs of Theorems}

\subsection{Proof of Theorem  \ref{Theorem: discrete EL equations 2}}  \label{Appendix: Discrete Proof 1}

\begin{theorem}  The Type~I discrete Hamilton's variational principle,
	\begin{equation*}
		\delta \bar{\mathfrak{S}}_d \left(\{  (q_k, \mathfrak{q}_k, \lambda _k) \}_{k=0}^{N} \right)  = 0,
	\end{equation*}
	where,
	\begin{equation*} 
		\bar{\mathfrak{S}}_d \left(\{  (q_k, \mathfrak{q}_k, \lambda _k) \}_{k=0}^{N} \right)  =  \sum_{k=0}^{N-1}{  \left[  L_d(q_k, \mathfrak{q}_k,  q_{k+1} , \mathfrak{q}_{k+1}) -  \lambda_{k} \frac{ \mathfrak{q}_{k+1} - \mathfrak{q}_k}{\tau_{k+1} - \tau_{k}}  + \lambda_k g(\mathfrak{q}_k)   \right]  \frac{ \mathfrak{q}_{k+1} - \mathfrak{q}_k}{\tau_{k+1} - \tau_{k}}    },
	\end{equation*}
	is equivalent to the discrete extended Euler--Lagrange equations,
	\small 
	\[   \mathfrak{q}_{k+1} =  \mathfrak{q}_k + (\tau_{k+1} - \tau_{k} )g(\mathfrak{q}_k)  ,\]   
	\[  \frac{ \mathfrak{q}_{k+1} - \mathfrak{q}_k}{\tau_{k+1} - \tau_{k}}   D_1 L_d(q_k, \mathfrak{q}_k,  q_{k+1} , \mathfrak{q}_{k+1}) + \frac{ \mathfrak{q}_{k} - \mathfrak{q}_{k-1}}{\tau_{k} - \tau_{k-1}}   D_3 L_d(q_{k-1}, \mathfrak{q}_{k-1},  q_{k} , \mathfrak{q}_{k})     = 0, \]
\begin{align*}
	&    \left[    D_2 L_{d_k} + \lambda_{k} \frac{ 1}{\tau_{k+1} - \tau_{k}}  + \lambda_k \nabla  g(\mathfrak{q}_k)  \right] \frac{ \mathfrak{q}_{k+1} - \mathfrak{q}_k}{\tau_{k+1} - \tau_{k}}      -  \frac{ 1}{\tau_{k+1} - \tau_{k}} \left[    L_{d_k} -  \lambda_{k} \frac{ \mathfrak{q}_{k+1} - \mathfrak{q}_k}{\tau_{k+1} - \tau_{k}}  + \lambda_k g(\mathfrak{q}_k)      \right]  \\ & \qquad    + \left[    D_4 L_{d_{k-1}} - \lambda_{k-1} \frac{ 1}{\tau_{k} - \tau_{k-1}}   \right] \frac{ \mathfrak{q}_{k} - \mathfrak{q}_{k-1}}{\tau_{k} - \tau_{k-1}}   +  \frac{ 1}{\tau_{k} - \tau_{k-1}} \left[     L_{d_{k-1}} -  \lambda_{k-1} \frac{ \mathfrak{q}_{k} - \mathfrak{q}_{k-1}}{\tau_{k} - \tau_{k-1}}  + \lambda_{k-1} g(\mathfrak{q}_{k-1})      \right]     =0,
\end{align*} 
where $L_{d_k} $ denotes $ L_d(q_k, \mathfrak{q}_k,  q_{k+1} , \mathfrak{q}_{k+1})$.
	\tiny 
	\proof{We use the notation $L_{d_k} = L_d(q_k, \mathfrak{q}_k,  q_{k+1} , \mathfrak{q}_{k+1})$, and we will use the fact that $\delta q_0 =\delta q_N =\delta \mathfrak{q}_0 =\delta \mathfrak{q}_N =0$ throughout the proof. We have
		\begin{align*}
			\delta \bar{\mathfrak{S}}_d & = \delta \left(  \sum_{k=0}^{N-1}{  \left[  L_d(q_k, \mathfrak{q}_k,  q_{k+1} , \mathfrak{q}_{k+1}) -  \lambda_{k} \frac{ \mathfrak{q}_{k+1} - \mathfrak{q}_k}{\tau_{k+1} - \tau_{k}}  + \lambda_k g(\mathfrak{q}_k)   \right]  \frac{ \mathfrak{q}_{k+1} - \mathfrak{q}_k}{\tau_{k+1} - \tau_{k}}    } \right)\\ &  = 			 \sum_{k=1}^{N-1}{ \left[    D_2 L_{d_k} + \lambda_{k} \frac{ 1}{\tau_{k+1} - \tau_{k}}  + \lambda_k \nabla  g(\mathfrak{q}_k)  \right] \frac{ \mathfrak{q}_{k+1} - \mathfrak{q}_k}{\tau_{k+1} - \tau_{k}}   \delta \mathfrak{q}_k }    - \sum_{k=1}^{N-1}{  \frac{ 1}{\tau_{k+1} - \tau_{k}} \left[     L_{d_k}  -  \lambda_{k} \frac{ \mathfrak{q}_{k+1} - \mathfrak{q}_k}{\tau_{k+1} - \tau_{k}}  + \lambda_k g(\mathfrak{q}_k)      \right]  \delta \mathfrak{q}_k }    \\ & \qquad \qquad  +\sum_{k=0}^{N-2}{ \left[    D_4 L_{d_k} - \lambda_{k} \frac{ 1}{\tau_{k+1} - \tau_{k}}   \right] \frac{ \mathfrak{q}_{k+1} - \mathfrak{q}_k}{\tau_{k+1} - \tau_{k}}   \delta \mathfrak{q}_{k+1} }   + \sum_{k=0}^{N-2}{  \frac{ 1}{\tau_{k+1} - \tau_{k}} \left[     L_{d_k} -  \lambda_{k} \frac{ \mathfrak{q}_{k+1} - \mathfrak{q}_k}{\tau_{k+1} - \tau_{k}}  + \lambda_k g(\mathfrak{q}_k)      \right]  \delta \mathfrak{q}_{k+1} }       
			\\ & \qquad \qquad    +   \sum_{k=1}^{N-1}{  \frac{ \mathfrak{q}_{k+1} - \mathfrak{q}_k}{\tau_{k+1} - \tau_{k}}   D_1 L_{d_k}     \delta q_{k} }    + \sum_{k=0}^{N-2}{   \frac{ \mathfrak{q}_{k+1} - \mathfrak{q}_k}{\tau_{k+1} - \tau_{k}}   D_3 L_{d_k}     \delta q_{k+1}     }   + \sum_{k=0}^{N-1}{ \frac{ \mathfrak{q}_{k+1} - \mathfrak{q}_k}{\tau_{k+1} - \tau_{k}}  \left(g(\mathfrak{q}_k) - \frac{ \mathfrak{q}_{k+1} - \mathfrak{q}_k}{\tau_{k+1} - \tau_{k}}  \right) \delta \lambda_k    } .
		\end{align*}
		Thus,
		\begin{align*}
			\delta \bar{\mathfrak{S}}_d &  = 			 \sum_{k=1}^{N-1}{ \left[    D_2 L_{d_k} + \lambda_{k} \frac{ 1}{\tau_{k+1} - \tau_{k}}  + \lambda_k \nabla  g(\mathfrak{q}_k)  \right] \frac{ \mathfrak{q}_{k+1} - \mathfrak{q}_k}{\tau_{k+1} - \tau_{k}}   \delta \mathfrak{q}_k }     - \sum_{k=1}^{N-1}{  \frac{ 1}{\tau_{k+1} - \tau_{k}} \left[     L_{d_k} -  \lambda_{k} \frac{ \mathfrak{q}_{k+1} - \mathfrak{q}_k}{\tau_{k+1} - \tau_{k}}  + \lambda_k g(\mathfrak{q}_k)      \right]  \delta \mathfrak{q}_k }    \\ & \qquad \qquad  +\sum_{k=1}^{N-1}{ \left[    D_4 L_{d_{k-1}}  - \lambda_{k-1} \frac{ 1}{\tau_{k} - \tau_{k-1}}   \right] \frac{ \mathfrak{q}_{k} - \mathfrak{q}_{k-1}}{\tau_{k} - \tau_{k-1}}   \delta \mathfrak{q}_{k} }     + \sum_{k=1}^{N-1}{  \frac{ 1}{\tau_{k} - \tau_{k-1}} \left[     L_{d_{k-1}} -  \lambda_{k-1} \frac{ \mathfrak{q}_{k} - \mathfrak{q}_{k-1}}{\tau_{k} - \tau_{k-1}}  + \lambda_{k-1} g(\mathfrak{q}_{k-1})      \right]  \delta \mathfrak{q}_{k} }       
			\\ & \qquad \qquad    +   \sum_{k=1}^{N-1}{ \left[  \frac{ \mathfrak{q}_{k+1} - \mathfrak{q}_k}{\tau_{k+1} - \tau_{k}}   D_1 L_{d_k} + \frac{ \mathfrak{q}_{k} - \mathfrak{q}_{k-1}}{\tau_{k} - \tau_{k-1}}   D_3 L_{d_{k-1}}  \right]      \delta q_{k} }   + \sum_{k=0}^{N-1}{ \frac{ \mathfrak{q}_{k+1} - \mathfrak{q}_k}{\tau_{k+1} - \tau_{k}}  \left(g(\mathfrak{q}_k) - \frac{ \mathfrak{q}_{k+1} - \mathfrak{q}_k}{\tau_{k+1} - \tau_{k}}  \right) \delta \lambda_k    } .
		\end{align*}
		As a consequence, if 
	\[   \mathfrak{q}_{k+1} =  \mathfrak{q}_k + (\tau_{k+1} - \tau_{k} )g(\mathfrak{q}_k)  ,\]   
	\[  \frac{ \mathfrak{q}_{k+1} - \mathfrak{q}_k}{\tau_{k+1} - \tau_{k}}   D_1 L_d(q_k, \mathfrak{q}_k,  q_{k+1} , \mathfrak{q}_{k+1}) + \frac{ \mathfrak{q}_{k} - \mathfrak{q}_{k-1}}{\tau_{k} - \tau_{k-1}}   D_3 L_d(q_{k-1}, \mathfrak{q}_{k-1},  q_{k} , \mathfrak{q}_{k})     = 0, \]
	\begin{align*}
		&    \left[    D_2 L_{d_k} + \lambda_{k} \frac{ 1}{\tau_{k+1} - \tau_{k}}  + \lambda_k \nabla  g(\mathfrak{q}_k)  \right] \frac{ \mathfrak{q}_{k+1} - \mathfrak{q}_k}{\tau_{k+1} - \tau_{k}}      -  \frac{ 1}{\tau_{k+1} - \tau_{k}} \left[    L_{d_k} -  \lambda_{k} \frac{ \mathfrak{q}_{k+1} - \mathfrak{q}_k}{\tau_{k+1} - \tau_{k}}  + \lambda_k g(\mathfrak{q}_k)      \right]  \\ & \qquad    + \left[    D_4 L_{d_{k-1}} - \lambda_{k-1} \frac{ 1}{\tau_{k} - \tau_{k-1}}   \right] \frac{ \mathfrak{q}_{k} - \mathfrak{q}_{k-1}}{\tau_{k} - \tau_{k-1}}   +  \frac{ 1}{\tau_{k} - \tau_{k-1}} \left[     L_{d_{k-1}} -  \lambda_{k-1} \frac{ \mathfrak{q}_{k} - \mathfrak{q}_{k-1}}{\tau_{k} - \tau_{k-1}}  + \lambda_{k-1} g(\mathfrak{q}_{k-1})      \right]     =0 ,
	\end{align*}
		then $	\delta \bar{\mathfrak{S}}_d \left(\{  (q_k, \mathfrak{q}_k, \lambda _k) \}_{k=0}^{N} \right)  = 0$. Conversely, if $	\delta \bar{\mathfrak{S}}_d \left(\{  (q_k, \mathfrak{q}_k, \lambda _k) \}_{k=0}^{N} \right)  = 0$, then a discrete fundamental theorem of the calculus of variations yields the above equations.
		\qed}
\end{theorem}

\subsection{Proof of Theorem  \ref{Theorem: discrete EL equations 1}}  \label{Appendix: Discrete Proof 2}

\begin{theorem}  The Type~I discrete Hamilton's variational principle,
	\begin{equation*}
		\delta \bar{\mathfrak{S}}_d \left(\{  (q_k, \mathfrak{q}_k, \lambda _k) \}_{k=0}^{N} \right)  = 0,
	\end{equation*}
	where,
	\begin{equation*} 
		\bar{\mathfrak{S}}_d \left(\{  (q_k, \mathfrak{q}_k, \lambda _k) \}_{k=0}^{N} \right)  =  \sum_{k=0}^{N-1}{\left\{ \frac{ \mathfrak{q}_{k+1} - \mathfrak{q}_k}{\tau_{k+1} - \tau_{k}}  \left[  L_d(q_k, \mathfrak{q}_k,  q_{k+1} , \mathfrak{q}_{k+1}) -  \lambda_{k} \right] + \lambda_k g(\mathfrak{q}_k)     \right\}},
	\end{equation*}
	is equivalent to the discrete extended Euler--Lagrange equations,
	\[   \mathfrak{q}_{k+1} =  \mathfrak{q}_k + (\tau_{k+1} - \tau_{k} )g(\mathfrak{q}_k),  \] 
	\[  \frac{ \mathfrak{q}_{k+1} - \mathfrak{q}_k}{\tau_{k+1} - \tau_{k}}  D_1 L_d(q_k, \mathfrak{q}_k,  q_{k+1} , \mathfrak{q}_{k+1}) +  \frac{ \mathfrak{q}_{k} - \mathfrak{q}_{k-1}}{\tau_{k} - \tau_{k-1}}  D_3 L_d(q_{k-1}, \mathfrak{q}_{k-1},  q_{k} , \mathfrak{q}_{k})  = 0, \]		\small 		\begin{align*}
	&  \frac{ \mathfrak{q}_{k+1} - \mathfrak{q}_k}{\tau_{k+1} - \tau_{k}}  D_2 L_{d_k} -  \frac{ 1}{\tau_{k+1} - \tau_{k}}L_{d_k}+ \frac{ \mathfrak{q}_{k} - \mathfrak{q}_{k-1}}{\tau_{k} - \tau_{k-1}}  D_4 L_{d_{k-1}} + \frac{ 1}{\tau_{k} - \tau_{k-1}}  L_{d_{k-1} } =  \frac{\lambda_{k-1}}{\tau_{k} - \tau_{k-1}}  - \frac{\lambda_k}{\tau_{k+1} - \tau_k}  - \lambda_k  \nabla  g(\mathfrak{q}_k)    ,
\end{align*}
where $L_{d_k} $ denotes $ L_d(q_k, \mathfrak{q}_k,  q_{k+1} , \mathfrak{q}_{k+1})$.
	\tiny	\proof{ We use the notation $L_{d_k} = L_d(q_k, \mathfrak{q}_k,  q_{k+1} , \mathfrak{q}_{k+1})$, and we will use the fact that $\delta q_0 =\delta q_N =\delta \mathfrak{q}_0 =\delta \mathfrak{q}_N =0$ throughout the proof. We have
		\begin{align*}
			\delta \bar{\mathfrak{S}}_d & = \delta \left(   \sum_{k=0}^{N-1}{\left\{ \frac{ \mathfrak{q}_{k+1} - \mathfrak{q}_k}{\tau_{k+1} - \tau_{k}}  \left[  L_d(q_k, \mathfrak{q}_k,  q_{k+1} , \mathfrak{q}_{k+1}) -  \lambda_{k} \right] + \lambda_k g(\mathfrak{q}_k)     \right\}} \right)\\ &  = 			 \sum_{k=1}^{N-1}{\left[ \frac{ \mathfrak{q}_{k+1} - \mathfrak{q}_k}{\tau_{k+1} - \tau_{k}}  D_2 L_{d_k}-  \frac{ 1}{\tau_{k+1} - \tau_{k}} L_{d_k}+ \frac{\lambda_k}{\tau_{k+1} - \tau_k}  + \lambda_k  \nabla g(\mathfrak{q}_k)   \right] \delta \mathfrak{q}_k }    + \sum_{k=0}^{N-2}{ \left[ \frac{ \mathfrak{q}_{k+1} - \mathfrak{q}_k}{\tau_{k+1} - \tau_{k}}  D_4 L_{d_k}  + \frac{ 1}{\tau_{k+1} - \tau_{k}}  L_{d_k} -  \frac{\lambda_k}{\tau_{k+1} - \tau_k}  \right] \delta \mathfrak{q}_{k+1}     }       
			\\ & \quad \qquad    +   \sum_{k=1}^{N-1}{   \frac{ \mathfrak{q}_{k+1} - \mathfrak{q}_k}{\tau_{k+1} - \tau_{k}}  D_1 L_{d_k} \delta q_{k} }    + \sum_{k=0}^{N-2}{    \frac{ \mathfrak{q}_{k+1} - \mathfrak{q}_k}{\tau_{k+1} - \tau_{k}}  D_3 L_{d_k} \delta q_{k+1}     }       + \sum_{k=0}^{N-1}{  \left(g(\mathfrak{q}_k) -  \frac{ \mathfrak{q}_{k+1} - \mathfrak{q}_k}{\tau_{k+1} - \tau_{k}} \right) \delta \lambda_k    } .
		\end{align*}
		Thus,
		\begin{align*}
			\delta \bar{\mathfrak{S}}_d & = 			 \sum_{k=1}^{N-1}{\left[ \frac{ \mathfrak{q}_{k+1} - \mathfrak{q}_k}{\tau_{k+1} - \tau_{k}}  D_2 L_{d_k} -  \frac{ 1}{\tau_{k+1} - \tau_{k}} L_{d_k} + \frac{\lambda_k}{\tau_{k+1} - \tau_k}  + \lambda_k  \nabla g(\mathfrak{q}_k) + \frac{ \mathfrak{q}_{k} - \mathfrak{q}_{k-1}}{\tau_{k} - \tau_{k-1}}  D_4 L_{d_{k-1}}  + \frac{ 1}{\tau_{k} - \tau_{k-1}}  L_{d_{k-1}} -  \frac{\lambda_{k-1}}{\tau_{k} - \tau_{k-1}}   \right] \delta \mathfrak{q}_k }       
			\\ & \qquad \qquad    +   \sum_{k=1}^{N-1}{\left[   \frac{ \mathfrak{q}_{k+1} - \mathfrak{q}_k}{\tau_{k+1} - \tau_{k}}  D_1 L_{d_k} +   \frac{ \mathfrak{q}_{k} - \mathfrak{q}_{k-1}}{\tau_{k} - \tau_{k-1}}  D_3 L_{d_{k-1}} \right]\delta q_{k} }        + \sum_{k=0}^{N-1}{  \left(g(\mathfrak{q}_k) -  \frac{ \mathfrak{q}_{k+1} - \mathfrak{q}_k}{\tau_{k+1} - \tau_{k}} \right) \delta \lambda_k    } . 
		\end{align*}
		As a consequence, if 
	\begin{align*}
		&  \frac{ \mathfrak{q}_{k+1} - \mathfrak{q}_k}{\tau_{k+1} - \tau_{k}}  D_2 L_{d_k} -  \frac{ 1}{\tau_{k+1} - \tau_{k}}L_{d_k}+ \frac{ \mathfrak{q}_{k} - \mathfrak{q}_{k-1}}{\tau_{k} - \tau_{k-1}}  D_4 L_{d_{k-1}} + \frac{ 1}{\tau_{k} - \tau_{k-1}}  L_{d_{k-1} } =  \frac{\lambda_{k-1}}{\tau_{k} - \tau_{k-1}}  - \frac{\lambda_k}{\tau_{k+1} - \tau_k}  - \lambda_k  \nabla  g(\mathfrak{q}_k)    ,
	\end{align*}
		\[  \frac{ \mathfrak{q}_{k+1} - \mathfrak{q}_k}{\tau_{k+1} - \tau_{k}}  D_1 L_d(q_k, \mathfrak{q}_k,  q_{k+1} , \mathfrak{q}_{k+1}) +  \frac{ \mathfrak{q}_{k} - \mathfrak{q}_{k-1}}{\tau_{k} - \tau_{k-1}}  D_3 L_d(q_{k-1}, \mathfrak{q}_{k-1},  q_{k} , \mathfrak{q}_{k})  = 0 ,\]
		\[   \mathfrak{q}_{k+1} =  \mathfrak{q}_k + (\tau_{k+1} - \tau_{k} )g(\mathfrak{q}_k),  \]   
		then $	\delta \bar{\mathfrak{S}}_d \left(\{  (q_k, \mathfrak{q}_k, \lambda _k) \}_{k=0}^{N} \right)  = 0$. Conversely, if $	\delta \bar{\mathfrak{S}}_d \left(\{  (q_k, \mathfrak{q}_k, \lambda _k) \}_{k=0}^{N} \right)  = 0$, then a discrete fundamental theorem of the calculus of variations yields the above equations.
		\qed}
\end{theorem}

\subsection{Proof of Theorem  \ref{Theorem: discrete EL equations Lie Group}}  \label{Appendix: Discrete Proof Lie Group}

	\begin{theorem} The Type~I discrete Hamilton's variational principle,
	\begin{equation*}
		\delta \bar{\mathfrak{S}}_d \left(\{  (q_k, \mathfrak{q}_k, \lambda _k) \}_{k=0}^{N} \right)  = 0,
	\end{equation*}
	where,
	\begin{equation*}
		\bar{\mathfrak{S}}_d \left(\{  (q_k, \mathfrak{q}_k, \lambda _k) \}_{k=0}^{N} \right) =   \sum_{k=0}^{N-1}{  \left[  L_d(q_k ,f_k , \mathfrak{q}_k, \mathfrak{q}_{k +1} )   -  \lambda_{k} \frac{ \mathfrak{q}_{k+1} - \mathfrak{q}_k}{\tau_{k+1} - \tau_{k}}  + \lambda_k g(\mathfrak{q}_k)   \right]  \frac{ \mathfrak{q}_{k+1} - \mathfrak{q}_k}{\tau_{k+1} - \tau_{k}},    }
	\end{equation*}
	is equivalent to the discrete extended Euler--Lagrange equations,
	\small 
	\[   \mathfrak{q}_{k+1} =  \mathfrak{q}_k + (\tau_{k+1} - \tau_{k} )g(\mathfrak{q}_k),  \]   
\begin{align*}   \emph{Ad}^*_{f_k^{-1}} \left( \emph{T}_e^* \emph{L}_{f_k}   D_2 L_{d_k}   \right) & =   \emph{T}_e^* \emph{L}_{q_k}  D_1  L_{d_k}     + \frac{\tau_{k+1} - \tau_{k}}{ \mathfrak{q}_{k+1} - \mathfrak{q}_k}     \frac{ \mathfrak{q}_{k} - \mathfrak{q}_{k-1}}{\tau_{k} - \tau_{k-1}}    \emph{T}_e^* \emph{L}_{f_{k-1}}   D_2 L_{d_{k-1}} ,       \end{align*}
\begin{align*}
	&    \left[    D_3 L_{d_k} + \lambda_{k} \frac{ 1}{\tau_{k+1} - \tau_{k}}  + \lambda_k \nabla  g(\mathfrak{q}_k)  \right] \frac{ \mathfrak{q}_{k+1} - \mathfrak{q}_k}{\tau_{k+1} - \tau_{k}}      -  \frac{ 1}{\tau_{k+1} - \tau_{k}} \left[    L_{d_k} -  \lambda_{k} \frac{ \mathfrak{q}_{k+1} - \mathfrak{q}_k}{\tau_{k+1} - \tau_{k}}  + \lambda_k g(\mathfrak{q}_k)      \right]  \\ & \qquad    + \left[    D_4L_{d_k}  - \lambda_{k-1} \frac{ 1}{\tau_{k} - \tau_{k-1}}   \right] \frac{ \mathfrak{q}_{k} - \mathfrak{q}_{k-1}}{\tau_{k} - \tau_{k-1}}   +  \frac{ 1}{\tau_{k} - \tau_{k-1}} \left[    L_{d_{k-1}}  -  \lambda_{k-1} \frac{ \mathfrak{q}_{k} - \mathfrak{q}_{k-1}}{\tau_{k} - \tau_{k-1}}  + \lambda_{k-1} g(\mathfrak{q}_{k-1})      \right]     =0,
\end{align*}
\normalsize
where $L_{d_k} $ denotes $ L_d(q_k ,f_k , \mathfrak{q}_k, \mathfrak{q}_{k +1} )   $. 
	\tiny 
	\proof{ We use the notation $L_{d_k} = L_d(q_k ,f_k , \mathfrak{q}_k, \mathfrak{q}_{k +1} )   $ and we will use the fact that $\delta q_0 = \delta q_N = \delta \mathfrak{q}_0 = \delta \mathfrak{q}_N =\eta_0 = \eta_N = 0$ throughout the proof. We have
		\begin{align*}
			\delta \bar{\mathfrak{S}}_d \left(\{  (q_k, \mathfrak{q}_k, \lambda _k) \}_{k=0}^{N} \right) & =   \delta \left( \sum_{k=0}^{N-1}{  \left[  L_d(q_k ,f_k , \mathfrak{q}_k, \mathfrak{q}_{k +1} )   -  \lambda_{k} \frac{ \mathfrak{q}_{k+1} - \mathfrak{q}_k}{\tau_{k+1} - \tau_{k}}  + \lambda_k g(\mathfrak{q}_k)   \right]  \frac{ \mathfrak{q}_{k+1} - \mathfrak{q}_k}{\tau_{k+1} - \tau_{k}}    } \right)  \\ & =   \sum_{k=1}^{N-1}{ \left[    D_3  L_{d_k} + \lambda_{k} \frac{ 1}{\tau_{k+1} - \tau_{k}}  + \lambda_k \nabla  g(\mathfrak{q}_k)  \right] \frac{ \mathfrak{q}_{k+1} - \mathfrak{q}_k}{\tau_{k+1} - \tau_{k}}   \delta \mathfrak{q}_k }    - \sum_{k=1}^{N-1}{  \frac{ 1}{\tau_{k+1} - \tau_{k}} \left[      L_{d_k}-  \lambda_{k} \frac{ \mathfrak{q}_{k+1} - \mathfrak{q}_k}{\tau_{k+1} - \tau_{k}}  + \lambda_k g(\mathfrak{q}_k)      \right]  \delta \mathfrak{q}_k }    \\ & \qquad \qquad  +\sum_{k=0}^{N-2}{ \left[    D_4  L_{d_k} - \lambda_{k} \frac{ 1}{\tau_{k+1} - \tau_{k}}   \right] \frac{ \mathfrak{q}_{k+1} - \mathfrak{q}_k}{\tau_{k+1} - \tau_{k}}   \delta \mathfrak{q}_{k+1} }     + \sum_{k=0}^{N-2}{  \frac{ 1}{\tau_{k+1} - \tau_{k}} \left[      L_{d_k}-  \lambda_{k} \frac{ \mathfrak{q}_{k+1} - \mathfrak{q}_k}{\tau_{k+1} - \tau_{k}}  + \lambda_k g(\mathfrak{q}_k)      \right]  \delta \mathfrak{q}_{k+1} }       
			\\ & \qquad \qquad    +   \sum_{k=1}^{N-1}{  \frac{ \mathfrak{q}_{k+1} - \mathfrak{q}_k}{\tau_{k+1} - \tau_{k}}   D_1  L_{d_k}   \delta q_{k} }  + \sum_{k=0}^{N-1}{   \frac{ \mathfrak{q}_{k+1} - \mathfrak{q}_k}{\tau_{k+1} - \tau_{k}}   D_2 L_{d_k} \delta f_{k}     }   + \sum_{k=0}^{N-1}{ \frac{ \mathfrak{q}_{k+1} - \mathfrak{q}_k}{\tau_{k+1} - \tau_{k}}  \left(g(\mathfrak{q}_k) - \frac{ \mathfrak{q}_{k+1} - \mathfrak{q}_k}{\tau_{k+1} - \tau_{k}}  \right) \delta \lambda_k    } .
		\end{align*}
	
		We can write $\delta g_k $ as $\delta g_k = g_k \eta_k$ for some $\eta_k \in \mathfrak{g}$. Then, taking the variation of the discrete kinematics equation $q_{k+1} = q_k f_k$ gives the equation $\delta q_{k+1} = \delta q_k f_k + q_k \delta f_k $ and $f_k = q_k^{-1} q_{k+1}$. Therefore,
		\begin{align*}
			\delta f_k & = q_k^{-1} \delta q_{k+1} - q_k^{-1} \delta q_k f_k   = q_k^{-1} q_{k+1} \eta_{k+1}  - q_k^{-1} q_{k} \eta_{k}f_k  = f_k \eta_{k+1} - \eta_k f_k,
		\end{align*} so
		\begin{align*}
			\delta \bar{\mathfrak{S}}_d \left(\{  (q_k, \mathfrak{q}_k, \lambda _k) \}_{k=0}^{N} \right) & =   \sum_{k=1}^{N-1}{ \left[    D_3 L_{d_k} + \lambda_{k} \frac{ 1}{\tau_{k+1} - \tau_{k}}  + \lambda_k \nabla  g(\mathfrak{q}_k)  \right] \frac{ \mathfrak{q}_{k+1} - \mathfrak{q}_k}{\tau_{k+1} - \tau_{k}}   \delta \mathfrak{q}_k }      - \sum_{k=1}^{N-1}{  \frac{ 1}{\tau_{k+1} - \tau_{k}} \left[      L_{d_k} -  \lambda_{k} \frac{ \mathfrak{q}_{k+1} - \mathfrak{q}_k}{\tau_{k+1} - \tau_{k}}  + \lambda_k g(\mathfrak{q}_k)      \right]  \delta \mathfrak{q}_k }    \\ & \qquad \quad +\sum_{k=1}^{N-1}{ \left[ \left(    D_4  L_{d_{k-1}} - \lambda_{k-1} \frac{ 1}{\tau_{k} - \tau_{k-1}}   \right) \frac{ \mathfrak{q}_{k} - \mathfrak{q}_{k-1}}{\tau_{k} - \tau_{k-1}}  +  \frac{ 1}{\tau_{k} - \tau_{k-1}} \left(      L_{d_{k-1}}-  \lambda_{k-1} \frac{ \mathfrak{q}_{k} - \mathfrak{q}_{k-1}}{\tau_{k} - \tau_{k-1}}  + \lambda_{k-1} g(\mathfrak{q}_{k-1})      \right)    \right]  \delta \mathfrak{q}_{k} }       
			\\ & \qquad \qquad  \quad +   \sum_{k=1}^{N-1}{  \frac{ \mathfrak{q}_{k+1} - \mathfrak{q}_k}{\tau_{k+1} - \tau_{k}}  \left( \text{T}_e^* \text{L}_{q_k}  D_1  L_{d_{k}}  \bullet  \eta_k \right)  }     + \sum_{k=0}^{N-1}{   \frac{ \mathfrak{q}_{k+1} - \mathfrak{q}_k}{\tau_{k+1} - \tau_{k}}   \left( \text{T}_e^* \text{L}_{f_k}   D_2L_{d_{k}}  \bullet  \left[   \eta_{k+1} - f_{k}^{-1} \eta_k f_k \right]   \right)    }    \\ & \qquad \qquad  \qquad + \sum_{k=0}^{N-1}{ \frac{ \mathfrak{q}_{k+1} - \mathfrak{q}_k}{\tau_{k+1} - \tau_{k}}  \left(g(\mathfrak{q}_k) - \frac{ \mathfrak{q}_{k+1} - \mathfrak{q}_k}{\tau_{k+1} - \tau_{k}}  \right) \delta \lambda_k    } .
		\end{align*}
		Then,
		\begin{align*}
			\delta \bar{\mathfrak{S}}_d \left(\{  (q_k, \mathfrak{q}_k, \lambda _k) \}_{k=0}^{N} \right) & =   \sum_{k=1}^{N-1}{ \left[    D_3  L_{d_k} + \lambda_{k} \frac{ 1}{\tau_{k+1} - \tau_{k}}  + \lambda_k \nabla  g(\mathfrak{q}_k)  \right] \frac{ \mathfrak{q}_{k+1} - \mathfrak{q}_k}{\tau_{k+1} - \tau_{k}}   \delta \mathfrak{q}_k }    - \sum_{k=1}^{N-1}{  \frac{ 1}{\tau_{k+1} - \tau_{k}} \left[      L_{d_k} -  \lambda_{k} \frac{ \mathfrak{q}_{k+1} - \mathfrak{q}_k}{\tau_{k+1} - \tau_{k}}  + \lambda_k g(\mathfrak{q}_k)      \right]  \delta \mathfrak{q}_k }    \\ & \qquad  +\sum_{k=1}^{N-1}{ \left[ \left(    D_4  L_{d_{k-1}} - \lambda_{k-1} \frac{ 1}{\tau_{k} - \tau_{k-1}}   \right) \frac{ \mathfrak{q}_{k} - \mathfrak{q}_{k-1}}{\tau_{k} - \tau_{k-1}}  +  \frac{ 1}{\tau_{k} - \tau_{k-1}} \left(      L_{d_{k-1}}-  \lambda_{k-1} \frac{ \mathfrak{q}_{k} - \mathfrak{q}_{k-1}}{\tau_{k} - \tau_{k-1}}  + \lambda_{k-1} g(\mathfrak{q}_{k-1})      \right)    \right]  \delta \mathfrak{q}_{k} }          
			\\ & \qquad \quad    +   \sum_{k=1}^{N-1}{  \frac{ \mathfrak{q}_{k+1} - \mathfrak{q}_k}{\tau_{k+1} - \tau_{k}}  \left( \text{T}_e^* \text{L}_{q_k}  D_1 L_{d_k}   \bullet  \eta_k \right)  }  + \sum_{k=0}^{N-1}{ \frac{ \mathfrak{q}_{k+1} - \mathfrak{q}_k}{\tau_{k+1} - \tau_{k}}  \left(g(\mathfrak{q}_k) - \frac{ \mathfrak{q}_{k+1} - \mathfrak{q}_k}{\tau_{k+1} - \tau_{k}}  \right) \delta \lambda_k    }   \\ & \qquad \qquad   + \sum_{k=0}^{N-1}{   \frac{ \mathfrak{q}_{k} - \mathfrak{q}_{k-1}}{\tau_{k} - \tau_{k-1}}   \left( \text{T}_e^* \text{L}_{f_{k-1}}   D_2 L_{d_{k-1}}   \bullet   \eta_{k}   \right)    }  - \sum_{k=0}^{N-1}{   \frac{ \mathfrak{q}_{k+1} - \mathfrak{q}_k}{\tau_{k+1} - \tau_{k}}   \left( \text{T}_e^* \text{L}_{f_k}   D_2 L_{d_k}  \bullet   \text{Ad}_{f_k^{-1}} \eta_k   \right)    }   .
		\end{align*}
		
		As a consequence, if 
			\[   \mathfrak{q}_{k+1} =  \mathfrak{q}_k + (\tau_{k+1} - \tau_{k} )g(\mathfrak{q}_k),  \]   
		\begin{align*}   \emph{Ad}^*_{f_k^{-1}} \left( \emph{T}_e^* \emph{L}_{f_k}   D_2 L_{d_k}   \right) & =   \emph{T}_e^* \emph{L}_{q_k}  D_1  L_{d_k}     + \frac{\tau_{k+1} - \tau_{k}}{ \mathfrak{q}_{k+1} - \mathfrak{q}_k}     \frac{ \mathfrak{q}_{k} - \mathfrak{q}_{k-1}}{\tau_{k} - \tau_{k-1}}    \emph{T}_e^* \emph{L}_{f_{k-1}}   D_2 L_{d_{k-1}}     ,   \end{align*}
		\begin{align*}
			&    \left[    D_3 L_{d_k} + \lambda_{k} \frac{ 1}{\tau_{k+1} - \tau_{k}}  + \lambda_k \nabla  g(\mathfrak{q}_k)  \right] \frac{ \mathfrak{q}_{k+1} - \mathfrak{q}_k}{\tau_{k+1} - \tau_{k}}      -  \frac{ 1}{\tau_{k+1} - \tau_{k}} \left[    L_{d_k} -  \lambda_{k} \frac{ \mathfrak{q}_{k+1} - \mathfrak{q}_k}{\tau_{k+1} - \tau_{k}}  + \lambda_k g(\mathfrak{q}_k)      \right]  \\ & \qquad    + \left[    D_4L_{d_k}  - \lambda_{k-1} \frac{ 1}{\tau_{k} - \tau_{k-1}}   \right] \frac{ \mathfrak{q}_{k} - \mathfrak{q}_{k-1}}{\tau_{k} - \tau_{k-1}}   +  \frac{ 1}{\tau_{k} - \tau_{k-1}} \left[    L_{d_{k-1}}  -  \lambda_{k-1} \frac{ \mathfrak{q}_{k} - \mathfrak{q}_{k-1}}{\tau_{k} - \tau_{k-1}}  + \lambda_{k-1} g(\mathfrak{q}_{k-1})      \right]     =0,
		\end{align*}
		then $	\delta \bar{\mathfrak{S}}_d \left(\{  (q_k, \mathfrak{q}_k, \lambda _k) \}_{k=0}^{N} \right)  = 0$. Conversely, if $	\delta \bar{\mathfrak{S}}_d \left(\{  (q_k, \mathfrak{q}_k, \lambda _k) \}_{k=0}^{N} \right)  = 0$, then a discrete fundamental theorem of the calculus of variations yields the above equations.
		\qed}
\end{theorem}

\hfill  \\

%\section*{Competing Interests}
The authors have no competing interests to declare that are relevant to the content of this article. \\

%\section*{Data Availability Statement}
The datasets generated during and/or analyzed during the current study are available from the corresponding author on reasonable request. \\

\hfill 

\bibliography{LagrangianAcceleration}
\bibliographystyle{plainnat}

\end{document}